\documentclass {article}

\usepackage{amsfonts,amsmath,latexsym}
\usepackage{amstext,amssymb,esint}
\usepackage{amsthm}
\usepackage {epsf}
\usepackage {epsfig}
\usepackage{tikz}
\usetikzlibrary{arrows, patterns}
\usepackage[scanall]{psfrag}
\usepackage {graphicx}
\usepackage[section]{placeins}
\usepackage{verbatim}
\usepackage{appendix}

\newcommand {\R}{\mathbf{R}}
\newcommand {\Ss}{\mathbf{S}}
\newcommand {\C}{\mathbf{C}}

\newcommand {\B}{\mathbf{B}}

\newcommand{\Ric}{\operatorname{Ric}}

\newtheorem {thm} {Theorem}
\newtheorem {prop} [thm] {Proposition}
\newtheorem {lemma} [thm] {Lemma}
\newtheorem {cor} [thm] {Corollary}
\newtheorem {defn} {Definition}
\newtheorem {rmk} {Remark}

\begin {document}

\title {On constant $Q$-curvature metrics with isolated singularities}
\author {Jesse Ratzkin
\footnote {Institut f\"ur Mathematik,
		Universit\"at W\"urzburg, {\tt jesse.ratzkin@mathematik.uni-wuerzburg.de}}}

\maketitle

\begin{abstract} In this paper we derive a refined asymptotic expansion, 
near an isolated singularity,  
for conformally flat metrics with constant positive $Q$-curvature and positive 
scalar curvature. The condition that the metric has constant $Q$-curvature 
forces the conformal factor to satisfy a fourth order nonlinear partial differential 
equation with critical Sobolev growth, whose leading term is the bilaplacian. 
We model our results on a similar asymptotic expansion for conformally 
flat, constant scalar curvature metrics proven by Korevaar, Mazzeo, Pacard, 
and Schoen. Along the way we analyze the linearization of the $Q$-curvature 
equation about the Delaunay metrics recently discovered by Frank and K\"onig, 
which may be of independent interest. 
\end{abstract}

\section{Introducion} 

Let $(M,g)$ be a Riemannian manifold of dimension $n \geq 5$. In this
paper we investigate the behavior of a conformal metric $\widetilde 
g = u^{\frac{4}{n-4}} g$ near an isolated singularity, subject to 
curvature conditions, namely constant and positive $Q$-curvature and 
positive scalar curvature. 

To begin we fix some notation. Let $R_g$ and $\Ric_g$ 
denote the scalar and Ricci curvature of $g$, and let 
$\Delta_g$ denote the Laplace-Beltrami operator. The 
(fourth-order) $Q$-curvature of $g$ is 
\begin{equation} \label{q_defn1}
Q_g = - \frac{1}{2(n-1)} \Delta_g R_g + \frac{n^3-4n^2+16n-16}
{8(n-1)^2(n-2)^2} R_g^2 - \frac{2}{(n-2)^2} |\Ric_g|^2. 
\end{equation} 
We can simplify this expression using the Schouten 
tensor
\begin{equation} \label{schouten_defn} 
A_g = \frac{1}{n-2} \left ( \Ric_g - \frac{R_g g}{2(n-1)} \right ), 
\qquad J_g = \operatorname{tr}_g (A_g) = \frac{R_g}{2(n-1)}, 
\end{equation} 
so that \eqref{q_defn1} becomes 
\begin{equation} \label{q_defn2} 
Q_g = -\Delta_g J_g - 2 |A_g |^2 + \frac{n}{2} J_g^2 . 
\end{equation} 
Associated to $Q_g$ we find the fourth order differential 
operator 
\begin{equation} \label{paneitz_op_defn} 
P_g (u) = (-\Delta_g)^2 u + \operatorname{div} 
\left ( 4 A_g (\nabla u, \cdot) - (n-2) J_g \nabla u \right ) + 
\frac{n-4}{2} Q_g, 
\end{equation} 
which enjoys the transformation rule 
\begin{equation} \label{conf_trans_rule1} 
\widetilde g = u^{\frac{4}{n-4}} g \Rightarrow 
P_{\widetilde g} (v) = u^{-\frac{n+4}{n-4}} P_g(uv). 
\end{equation} 
Substituting $v=1$ into \eqref{conf_trans_rule1} 
we find 
\begin{equation} \label{conf_trans_rule2} 
\widetilde g = u^{\frac{4}{n-4}} \Rightarrow 
Q_{\widetilde g} = \frac{2}{n-4} u^{-\frac{n+4}{n-4}} 
P_g(u). 
\end{equation}  

S. Paneitz \cite{Pan1} first introduced the 
operator \eqref{paneitz_op_defn} and studied its transformation 
properties. Later T. Branson \cite{Bran1, Bran2} extended this 
operator to differential forms and studied a related 
sixth-order operator, as well as the $Q$-curvature. The 
reader can find summaries of the current understanding of $Q$-curvature in 
the  survey articles \cite{BG}, \cite{CEOY}, and \cite{HY}. 

\subsection{An aside on scalar curvature and the Yamabe problem} 

A this point we pause to discuss the Yamabe problem, which provides 
a guide for much of the investigation of the properties of 
$Q$-curvature. One can define the conformal Laplacian 
\begin{equation} \label{conf_lap}
\mathcal{L}_g (u) = -\Delta_g (u) + \frac{4(n-1)}{n-2} R_g u, 
\end{equation}
which enjoys the transformation rule 
\begin{equation} \label{conf_trans_rule3} 
\widetilde g = u^{\frac{4}{n-2}} g \Rightarrow 
\mathcal{L}_{\widetilde g} (v) = u^{-\frac{n+2}{n-2}} 
\mathcal{L}_g(uv).
\end{equation} 
Substituting $v=1$ into \eqref{conf_trans_rule3} 
we obtain 
\begin{equation} \label{conf_trans_rule4} 
\widetilde g = u^{\frac{4}{n-2}} g \Rightarrow 
R_{\widetilde g} = \frac{n-2}{4(n-1)} u^{-\frac{n+2}{n-2}}
\mathcal{L}_g(u).
\end{equation} 

This last transformation rule allows us to define 
the conformal invariant 
\begin{eqnarray} \label{yamabe_inv1} 
\mathcal{Y} ([g]) & = & \inf \left \{  \frac{ \int_M 
R_{\widetilde g} d\mu_{\widetilde g}} 
{(\operatorname{Vol}_{\widetilde g}(M))^{\frac{n-4}{n}}} : \widetilde g = 
u^{\frac{4}{n-2}} g \in [g] \right \} \\ \nonumber 
& =& \inf \left \{ \frac{n-2}{4(n-1)}\frac{ \int_M u \mathcal{L}_g(u)  
d\mu_g} {\left ( \int_M u^{\frac{2n}{n-2}} d\mu_g \right )^{\frac{n-2}{n}}} : 
u \in \mathcal{C}^\infty (M), u>0 \right \}
\end{eqnarray}
A straightforward computation shows that critical points of the 
functional 
$$\widetilde g\in [g] \mapsto \frac{\int_M R_{\widetilde g} 
d\mu_{\widetilde g}}{(\operatorname{Vol}_{\widetilde g}
(M))^{\frac{n-4}{n}}} $$
are precisely the constant scalar curvature metrics in 
the conformal class $[g]$.

In \cite{Y} Yamabe proposed finding a constant scalar curvature metric 
in a given conformal class $[g]$ through the family of related variational 
problems 
\begin{equation} \label{yamabe_inv_p} 
\mathcal{S}_p ([g]) = 
\inf \left \{ \frac{n-2}{4(n-1)} \frac{\int_M u \mathcal{L}_g(u)
d\mu_g}{\left ( \int_M u^p d\mu_g \right )^{2/p}} : 
u \in \mathcal{C}^\infty(M), u>0  \right \}
\end{equation} 
for $p$ such that $1 < p < \frac{2n}{n-2}$. By Rellich's compactness 
theorem the infimum $\mathcal{S}_p ([g])$ is a minimum, and 
realized by a smooth metric. Yamabe's strategy was to find a 
constant scalar curvature metric as a limit of minimizers 
for $\mathcal{S}_p([g])$ as $p \nearrow \frac{2n}{n-2}$. 

Trudinger \cite{Tru} first pointed out the difficulties in extracting
this limit, primarily that Rellich's compactness fails exactly 
for $p = \frac{2n}{n-2}$. Later Aubin \cite{Aub} resolved 
these difficulties in many cases, and Schoen \cite{Sch} 
completed Yamabe's original program. Schoen's resolution of 
the Yamabe problem culminated the work of many people over 
25 years, and continues to inspire new research. 

One can see the lack of compactness 
directly in the case of the round metric $g_0$ on the 
sphere $\Ss^n$. Using stereographic projection we can write 
the round metric as 
\begin{equation} \label{scal_curv_sph1}
g_0 = \frac{4}{(1+|x|^2)^2} \delta = \left ( \left ( \frac{1+|x|^2}{2}
\right )^{\frac{2-n}{2}}\right )^{\frac{4}{n-2}} \delta = 
u_1^{\frac{4}{n-2}} \delta
\end{equation} 
where $\delta$ is the flat metric. Dilating the 
Euclidean coordinates by a factor of $\lambda>0$ one finds the 
conformal factor 
\begin{equation} \label{scal_curv_sph2}
u_\lambda =\left ( \frac{1+\lambda^2 |x|^2}{2\lambda} 
\right )^{\frac{2-n}{2}}.
\end{equation}
However, as $\lambda \rightarrow \infty $ we see that $u_\lambda 
\rightarrow 0$ outside of any fixed neightborhood of $0$, 
while $u_\lambda(0) \rightarrow \infty$. Geometrically, this 
blow-up behavior concentrates the entirety of the sphere in 
a small neighborhood of the south pole, which 
corresponds to the origin in Euclidean coordinates, shrinking 
the complement of this neighborhood to be vanishingly small. 

The blow-up described above motivates one to understand 
the asymptotics of constant scalar curvature metrics with 
isolated singularities. In the special case that the 
background metric $g$ is locally conformally flat one can 
(locally) write a constant scalar curvature metric in the conformal 
class $[g]$ as $u^{\frac{4}{n-2}} \delta$. A theorem 
of Schoen and Yau \cite{SY} implies the value of the scalar 
curvature must be a positive constant, which we 
normalize to be $n(n-1)$. In this case  case 
\eqref{conf_trans_rule4} becomes 
\begin{equation} \label{scal_curv_pde}
-\Delta u = \frac{n(n-2)}{4} u^{\frac{n+2}{n-2}} .
\end{equation} 
A computation shows that $u_\lambda$, defined in \eqref{scal_curv_sph2}
solves \eqref{scal_curv_pde}, and indicates that one should 
study the possible asymptotic behavior of positive solutions of 
\eqref{scal_curv_pde}. In \cite{CGS} Caffarelli, Gidas and Spruck 
proved that any positive solution of \eqref{scal_curv_pde} in the 
punctured ball $\B_1(0) \backslash \{ 0 \}$ satisfies 
$$u(x) = \bar{u}(|x|) (1+ \mathcal{O}(|x|)), \qquad 
\bar{u} (r) = \frac{1}{r^{n-1} |\Ss^{n-1}|} \int_{|x| = r} 
u(x) d\sigma(x).$$
Later Korevaar, Mazzeo, Pacard and Schoen \cite{KMPS} 
refined this asymptotic expansion, obtaining the next 
term in the expansion and giving the terms in this expansion a 
more geometric interpretation. 

\subsection{Motivation for studying constant $Q$-curvature 
metrics with isolated singularities}

Returning to our discussion of $Q$-curvature, we can 
follow the well-established model of scalar curvature 
and define the conformal invariant 
\begin{eqnarray} \label{paneitz_inv1} 
\mathcal{Y}_4^+ ([g])& = & \inf\left \{ \frac{\int_M 
Q_{\widetilde g} d\mu_g}{(\operatorname{Vol}_{\widetilde 
g} (M))^{\frac{n-4}{n}} } : \widetilde g = u^{\frac{4}{n-4}} g \in 
[g] \right \} \\ \nonumber 
& = & \inf \left \{ \frac{2}{n-4} \frac{\int_M u P_g (u) 
d\mu_g}{ \left ( \int_M u^{\frac{2n}{n-4}} d\mu_g
\right )^{\frac{n-4}{n}} } : u \in \mathcal{C}^\infty
(M) , u> 0\right \}. 
\end{eqnarray} 
Once again, a straight-forward computation implies critical points 
of the functional 
$$u \mapsto \frac{\int_M u P_g (u) d\mu_g}{\left ( \int_M 
u^{\frac{2n}{n-4}} d\mu_g \right )^{\frac{n-4}{n}}} $$ 
give exactly the constant $Q$-curvature metrics within 
the conformal class $[g]$.  In the case that this constant is positive, 
we normalize it to be $\frac{n(n^2-4)}{8}$, which is the value attained 
by the round metric on the sphere. 

In constrast to the situation with 
scalar curvature, the existence and properties of minimizers and 
higher order critical points is poorly understood. However, 
the conformal class of the round metric on the sphere still 
provides an illustrative example. We once again write the 
round metric as a conformally flat metric using stereographic 
projection, so that 
\begin{equation} \label{sph_soln1} 
g_0 = \frac{4}{(1+|x|^2)^2} \delta = \left ( \left ( \frac{1+|x|^2}
{2} \right )^{\frac{4-n}{2}} \right )^{\frac{4}{n-4}} \delta = 
U_1(x)^{\frac{4}{n-4}}\delta.
\end{equation}
Once again we may apply a conformal dilation, which gives the 
solution 
\begin{equation} \label{sph_soln2} 
U_\lambda = \left ( \frac{1+\lambda^2 |x|^2}{2\lambda} 
\right )^{\frac{4-n}{2}}
\end{equation} 
for any $\lambda >0$. 

We can seek metrics $g = u^{\frac{4}{n-4}} g_0$ in the 
conformal class of the round metric $g_0$ with constant 
$Q$-curvature $\frac{n(n^2-4)}{8}$. In steregraphic 
coordinates \eqref{conf_trans_rule2} reduces to 
\begin{equation} \label{paneitz_pde1} 
\Delta^2 u = \frac{n(n-4)(n^2-4)}{16} u^{\frac{n+4}{n-4}}.
\end{equation} 
By the above discussion $U_\lambda$ solves \eqref{paneitz_pde1}
for each $\lambda > 0$. In fact, a theorem of C. S. Lin \cite{Lin}
states that any positive solution of \eqref{paneitz_pde1} must have 
the form $U_\lambda (x-x_0)$ for some $\lambda > 0$ and 
$x_0 \in \R^n$. 

As before, we examine the behavior of $U_{\lambda, x_0} = 
U_\lambda (\cdot + x_0)$ as $\lambda \rightarrow \infty$, finding 
that $U_{\lambda, x_0} \rightarrow 0$ outside any fixed neighborhood
of $x_0$, while it blows up at $x_0$. This behavior illustrates the 
natural lack of compactness in the infimum \eqref{paneitz_inv1}, 
and explains why we should study solutions of \eqref{paneitz_pde1} 
with isolated singularities. 

\subsection{Main results} 

We concentrate on solututions of \eqref{paneitz_pde1} in the 
unit ball $\B_1(0)$ with an isolated singularity at the origin.  

Recently Jin and Xiong \cite{JX} proved that if $u \in \mathcal{C}^\infty
(\B_1 (0) \backslash \{ 0 \})$ is a positive solution of \eqref{paneitz_pde1} 
which also satisfies $-\Delta u >0$ then 
\begin{equation} \label{JX_asymp}
u(x) = \bar u(|x|) (1+ \mathcal{O} (x)), \qquad \bar u (r) = 
\frac{1}{r^{n-1} |\Ss^{n-1}|} \int_{|x| = r} u(x) d\sigma (x).
\end{equation} 

To understand our results we introduce some 
special solutions of \eqref{paneitz_pde1}. We change to 
cylindrical coordinates, letting $t = -\log |x|$, $\theta = 
x/|x|$, and 
\begin{equation} \label{cyl_coord_change} 
v: (0,\infty) \times \Ss^{n-1} \rightarrow (0,\infty), \qquad 
v(t,\theta) = e^{\left ( \frac{4-n}{2}\right )t} u(e^{-t} \theta) . 
\end{equation} 
Under this change of coordinates \eqref{paneitz_pde1} becomes 
\begin{eqnarray} \label{paneitz_pde2} 
\frac{n(n-4)(n^2-4)}{16} v^{\frac{n+4}{n-4}} & = & \frac{\partial^4 v}{\partial t^4} 
- \left ( \frac{n(n-4)+8}{2} \right )\frac{\partial^2 v}{\partial t^2} + 
\frac{n^2 (n-4)^2}{16} v \\ \nonumber 
&& + \Delta_\theta^2 v + 2 \Delta_\theta \frac{\partial^2 v}{\partial t^2} 
- \frac{n(n-4)}{2}\Delta_\theta v, 
\end{eqnarray} 
where $\Delta_\theta$ is the Laplace-Beltrami operator on the round 
sphere $\Ss^{n-1}$, and the condition $-\Delta u >0$ becomes 
\begin{equation} \label{pos_lap_cyl_coords} 
-\ddot v + 2\dot v + \frac{n(n-4)}{4} v - \Delta_\theta v > 0.
\end{equation} 

In \cite{Lin}, Lin also proved that any global solution of \eqref{paneitz_pde2} 
must be a function of $t$ alone, and so it must satisfy the ordinary differential 
equation (ODE) 
\begin{equation} \label{paneitz_ode1} 
\frac{n(n-4)(n^2-4)}{16} v^{\frac{n+4}{n-4}} = \ddddot v - \left ( 
\frac{n(n-4)+8}{2} \right ) \ddot v + \frac{n^2(n-4)^2}{16} v . 
\end{equation} 
One immediately sees two solutions of \eqref{paneitz_ode1}: 
the constant solution 
\begin{equation} \label{cyl_soln}
v_{cyl} = \left ( \frac{n(n-4)}{n^2-4} \right )^{\frac{n-4}{8}}
\end{equation} 
and the spherical solution 
\begin{equation} 
u(x) = U_1(x) = \left ( \frac{1+|x|^2}{2} \right )^{\frac{4-n}{2}}
\Leftrightarrow v_{sph} = (\cosh t)^{\frac{4-n}{2}} . 
\end{equation} 

Frank and K\"onig \cite{FK} classified all global, positive solutions of 
\eqref{paneitz_ode1}. First they demonstrate the existence of a unique  
periodic solution $v_\epsilon$ of \eqref{paneitz_ode1} attaining 
its minimum value of $\epsilon$ at $t=0$ for each $\epsilon \in 
(0,v_{cyl}]$. We call $v_\epsilon$ the {\bf Delaunay solution} 
with necksize $\epsilon$. In the same paper they prove that 
each global, positive 
solution of \eqref{paneitz_ode1} must either have the form 
$(\cosh (\cdot + T))^{\frac{4-n}{2}}$ or have the form $v_\epsilon 
(\cdot + T)$ for some $\epsilon \in (0, v_{cyl}]$ and $T \in \R$. 

We are now ready to state our main theorem. 
\begin{thm} \label{refined_asymp_thm} 
Let $v\in \mathcal{C}^\infty((0,\infty) \times \Ss^{n-1})$ be a 
positive solution of \eqref{paneitz_pde2} which also 
satisfies \eqref{pos_lap_cyl_coords}. Then 
either $\limsup_{t \rightarrow \infty} v(t,\theta) = 0$ or 
there exists parameters $\epsilon \in (0,v_{cyl}]$, $T \in \R$ 
and $a \in \R^n$, 
and positive constants $C$ and $\beta > 1$ such that 
\begin{equation} 
\left | v(t,\theta) - v_\epsilon (t+T)- e^{-t} \langle \theta, a \rangle 
\left ( -\dot v_\epsilon (t+T) + \frac{n-4}{2} v_\epsilon (t+T) \right ) 
\right | \leq C e^{-\beta t} . 
\end{equation} 
\end{thm} 
\begin{rmk} We will see below in \eqref{trans_del_soln2} and 
\eqref{trans_del_soln3} that one obtains 
$$v_\epsilon + e^{-t} \langle \theta, a \rangle \left ( -\dot v_\epsilon 
+ \frac{n-4}{2} v_\epsilon \right )$$ 
by applying a certain translation to the Delaunay solution $v_\epsilon$. 
\end{rmk}

\begin{rmk} \label{further_asymp_rmk}
We speculate that one can adapt our 
proof below to prove a further refinement of our 
estimate in Theorem \ref{refined_asymp_thm}, following 
the estimate of Han, Li and Li \cite{HLL}. In \eqref{indicial_defn} we 
define the increasing set of indicial roots 
$$\Gamma_\epsilon = \{ \dots, -\gamma_{\epsilon,2}, 
-\gamma_{\epsilon,1} = -1, 0, \gamma_{\epsilon,1} = 1, 
\gamma_{\epsilon,2}, \dots \}, \qquad \gamma_{\epsilon,j} 
\rightarrow \infty.$$ 
One should be able to prove an estimate of the form 
$$
\left | v(t,\theta) - v_\epsilon (t+T) - \sum_{i=1}^m \sum_{j=0}^{m-1} 
c_{ij} (t,\theta) t^j e^{-\gamma_{\epsilon, j} t} \right | \leq 
C t^m e^{-\gamma_{\epsilon, m+1}}$$ 
where the coefficient functions $c_{ij} (t,\theta)$ are bounded. 
\end{rmk} 

We recast Theorem \ref{refined_asymp_thm} in geometric 
terms. 
\begin{cor} 
Let $n \geq 5$ and let $g = u^{\frac{4}{n-4}} \delta$ be a conformally 
flat metric on $\B_1(0) \backslash \{ 0 \}$ with positive scalar curvature 
and $Q_g = \frac{n(n^2-4)}{8}$. Then either $g$ extends to a smooth 
metric on $\B_1(0)$ or there exist parameters $\epsilon \in (0,v_{cyl}]$, 
$T \in \R$ and $a \in \R$ and $\beta >1$ such that 
$$u(x) = |x|^{\frac{4-n}{2}} \left ( v_\epsilon (-\log |x| + T) 
+ \langle x,a \rangle \left ( -\dot v_\epsilon (-\log |x|+T) + \frac{n-4}{2}
v_\epsilon (-\log |x| + T) \right )  + \mathcal{O} (|x|^\beta) \right ) $$
\end{cor} 

\begin{proof} We have already seen that $Q_g = \frac{n(n^2-4)}{8}$ 
is equivalent to \eqref{paneitz_pde1}. By \eqref{conf_trans_rule4} 
$$R_g = \frac{n-2}{4(n-1)}u^{-\left ( \frac{n+2}{n-4} \right )}
\left ( -\Delta u^{\frac{n-2}{n-4}} \right ) ,$$
so that 
\begin{equation} \label{pos_scal_curv}
R_g > 0 \Leftrightarrow -\Delta u^{\frac{n-2}{n-4}} > 0 \Leftrightarrow 
-\Delta u > \frac{2}{n-4} \frac{|\nabla u|^2}{u}. 
\end{equation} 
In particular, $R_g >0$ implies $-\Delta u >0$. Changing to cylindrical 
coordinates we obtain a function 
$$v \in \mathcal{C}^\infty ((0,\infty) \times \Ss^{n-1}), \qquad 
v(t,\theta) = e^{\left ( \frac{4-n}{2} \right ) t} u(e^{-t} \theta)$$ 
which satisfies \eqref{paneitz_pde2} and \eqref{pos_lap_cyl_coords}.
\end{proof}

We constrast our methods with that of Jin and Xiong \cite{JX}. 
They use Green's identity to transform \eqref{paneitz_pde1} 
into an integral formula
\begin{equation} \label{paneitz_int_form1} 
u(x) = \int_{\B_1(0)} \frac{(u(y))^{\frac{n+4}{n-4}}}{|x-y|^{n-4}}
d\mu(y) + h(x) , \qquad u \in L^{\frac{n+4}{n-2}} (\B_1(0)) 
\cap \mathcal{C}^1(\B_1(0) \backslash \{ 0 \})
\end{equation} 
where $h \in \mathcal{C}^1(\B_1(0))$. In this setting the condition 
$-\Delta u > 0$ implies $h>0$, which allows Jin and Xiong to apply 
the method of moving spheres/planes. They first prove {\it a priori} 
upper and lower bounds for $u$, and then prove their asymptoic 
estimate using moving spheres. 

Our proof below follows the techniques in \cite{KMPS}, extracting 
a limit from a slide-back sequence. We start with $v:(0,\infty) \times \Ss^{n-1} 
\rightarrow (0,\infty)$ which satisfies \eqref{paneitz_pde2} and a 
sequence $\tau_k \rightarrow \infty$. and let $v_k (t,\theta) = 
v(t+\tau_k, \theta)$. The {\it a priori} estimates of Jin and Xiong 
allow us to extract a convergent subsequence, which is defined 
on all of $\R \times \Ss^{n-1}$. By the uniqueness theorem of 
Frank and K\"onig this limit must have the form $v_\epsilon (t + T)$
for some $\epsilon \in (0, v_{cyl}]$ and $T \in \R$. The main task 
in proving simple asymptotics involves proving $\epsilon$ and 
$T$ do not depend on the choice of the sequence $\tau_k 
\rightarrow \infty$ or the choice of the convergent subsequence. 
Our proof that the parameters $\epsilon$ and $T$ are independent 
of all choices relies heavily on a careful analysis of the linearization 
of the PDE \eqref{paneitz_pde1} about a Delaunay solution, including 
an asymptotic expansion of solutions of the linearized equation. Once 
we prove the simple asymptotics, we obtain the refined asymptotics 
from the asymptotic expansion of the linearized operator. 

The two techniques are complementary. The proof of Jin and Xiong 
is more general, and applies to any solution of the integral 
equation 
$$u(x) = \int_{\B_1(0)} \frac{(u(y))^{\frac{n+2\sigma}{n-2\sigma}}}
{|x-y|^{n-2\sigma}} d\mu(y) + h(x), \qquad h>0, \quad 0 < \sigma
< \frac{n}{2}, $$
whereas at this time one can apply our technique only in the 
cases $\sigma=1$ and $\sigma =2$. On the other hand, their technique 
does not give the refined asymptotic expansion. 

Previously Gonz\'alez \cite{Gonz} proved a similar asymptotics 
theorem for $\sigma_k$-curvature, which is a fully-nonlinear 
genearlization of scalar curvature whose associated PDE is second order and 
fully nonlinear. More recently, Caffarelli, Jin, Sire, and Xiong \cite{CJSX}
prove an asymptotic result for positive solutions of the nonlocal equation 
$$(-\Delta)^s u = u^{\frac{n+2s}{n-2s}}$$
for any $s \in (0,1)$. Also Baraket and Rebhi \cite{BR} and 
Y.-J. Lin \cite{YJLin} construct many examples of constant $Q$-curvature 
metrics. 

The rest of this paper proceeds as follows. Section \ref{prelim_sec} we 
collect some preliminary computations, most of which exist already 
in the literature. Within this section we include a detailed description 
of the Delaunay solutions. Our analysis begins in earnest in 
Section \ref{lin_anal_sec} where we study the linearization of 
the PDE \eqref{paneitz_pde2}. Most importantly we prove the linear 
stability of the Delaunay solutions using the Fourier-Laplace transform 
as developed by Mazzeo, Pollack and Uhlenbeck \cite{MPU}. In this 
section we also define the indicial roots mentioned above in 
Remark \ref{further_asymp_rmk}, which give the exponential 
growth rates of the solutions of the linearization of \eqref{paneitz_pde2} 
when linearized about a Delaunay solution. In Section \ref{asymp_sec} 
we present an alternative proof of the simple asymptotics Jin and 
Xiong prove in \cite{JX}, and in Section \ref{refined_sec} we 
derive the refined asymptotics. 

{\sc Acknowledgements:} This research was partially supported 
by the National Research Foundation of South Africa. I would like to 
thank Rupert Frank, Andrea Malchiodi, and Rafe Mazzeo for 
helpful conversations during the conference ``Recent Advances in Nonlocal 
and Nonlinear Analysis" at ETH Z\"urich in June 2014, as well as the organizers 
of this conference for providing a fruitful venue. 

\section{Preliminaries} \label{prelim_sec} 

We collect some preliminary computations related to solutions 
of \eqref{paneitz_pde1} and constant $Q$-curvature metrics. 

\subsection{Symmetries of the PDE}

In this section we recall some of the symmetries of \eqref{paneitz_pde1}, 
which essentially arise from the transformation rule \eqref{conf_trans_rule1}.
One can of course translate solutions to obtain a new solution, but two 
more interesting symmetries reflect the scale invariance and conformal 
invariance outlined above. 

We first discuss scale invariance. Let $u>0$ solve \eqref{paneitz_pde1}
and let $\lambda>0$. We seek $a>0$ such that $u_\lambda(x) = 
\lambda^a u(\lambda x)$ is also a solution. Evaluating, we find 
$$(-\Delta)^2 u_\lambda = \lambda^{a+4} (-\Delta)^2 u (\lambda x), \qquad 
\frac{n(n-4)(n^2-4)}{16} u_\lambda^{\frac{n+4}{n-4}} = \lambda^{\frac{(n+4)a}{n-4}}
\frac{n(n-4)(n^2-4)}{16} (u(\lambda x))^{\frac{n+4}{n-4}},$$
which coincide precisely when $a = \frac{n-4}{2}$. In other words we see 
\begin{equation} \label{scaling_law1}
u \textrm{ solves \eqref{paneitz_pde1} } \Rightarrow 
u_\lambda(x) = \lambda^{\frac{n-4}{2}} u(\lambda x) 
\textrm{ solves \eqref{paneitz_pde1} for all } \lambda >0.
\end{equation} 

This symmetry is much simpler in the cylindrical coordinates. Letting 
$T = -\log \lambda$ we see 
$$
v_\lambda (t,\theta) = e^{\frac{(4-n)t}{2}} u_\lambda (e^{-t} \theta) =
\lambda^{\frac{n-4}{2}} e^{\frac{(4-n)t}{2}} u(\lambda e^{-t} \theta) 
= e^{\frac{(4-n)(t+T)}{2}} u(e^{-(t+T)} \theta),$$
so that 
$$
v \textrm{ solves \eqref{paneitz_pde2} } \Rightarrow 
v_T(t,\theta) = v(t+T,\theta) 
\textrm{ solves \eqref{paneitz_pde2} for all } T\in \R,
$$
which is readily apparent directly from the PDE. 

The second symmetry reflects the invariance under reflections 
through spheres, and one can write 
write it explicitly by defining the Kelvin transforms
\begin {equation} \label {kelvin_trans}
\mathbb{K}_{x_0} (u) (x) = |x|^{2-n} u \left ( \frac{x}{|x|^2} + x_0  
\right ), \qquad \widehat{\mathbb{K}}_{x_0} (u)(x) = |x|^{4-n} u \left 
( \frac{x}{|x|^2} + x_0 \right ) . \end {equation}
Observe that the functions $\mathbb{K}_{x_0}(u)$ and 
$\widehat {\mathbb{K}}_{x_0}
(u)$ are now defined on different domains. For instance, if $u \in 
\mathcal{C}^\infty (\B_1(0) \backslash \{ 0 \} )$ then 
$\mathbb{K}_{x_0} (u), \widehat {\mathbb{K}}_{x_0} \in 
\mathcal{C}^\infty (\R^n \backslash \B_1 (x_0))$.

One can find the usual Kelvin transformation law 
\begin {equation} \label{orig_kelvin_law} 
\Delta (\mathbb{K}_{x_0}(u))(x) = |x|^{-4} \mathbb{K}_0(\Delta u)(x)
\end {equation}
in many textbooks, and the transformation law 
\begin {equation} \label{4th_order_kelvin}
\Delta^2 (\widehat{\mathbb{K}}_{x_0}(u))(x) = |x|^{-8} \widehat{\mathbb{K}}_{x_0}
(\Delta^2 u)(x), \end {equation}
appears (for instance) in Lemma 3.6 of \cite{Xu1}. 

\begin{rmk} In the case $x_0 =0$ 
this transformation looks particularly simple in cylindrical coordinates, 
namely $v(t,\theta)$ gets transformed to $\widehat{\mathbb{K}}_0(v)(t,\theta)
= v(-t,\theta)$. On the other hand, the transformation is much more 
complicated in cylindrical coordinates when the center is not $0$.
\end{rmk}

\subsection{Delaunay solutions} 

We have already introduced the Delaunay solutions, and 
in this section we give a more detailed description. 

The family of Delaunay solutions account for all 
of the positive solutions of \eqref{paneitz_pde2} on 
the whole cylinder $\R \times \Ss^{n-1}$. Recall the 
cylindrical (constant) and spherical solutions 
$$ v_{cyl} = \left ( \frac{n(n-4)}{n^2-4} \right )^{\frac{n-4}{8}}, 
\qquad v_{sph} = (\cosh t)^{\frac{4-n}{2}}.$$
It is convenient to observe that, because $n > 4$ we have 
$$n^2-4 -n(n-4) = 4n-4 >0 \Rightarrow 0<v_{cyl} < 1,$$
as well as 
$$v_{sph} (0) = 1, \qquad \dot v_{sph} (t) > 0 \textrm{ for }
t<0, \qquad \dot v_{sph}(t) < 0 \textrm{ for }t>0, \qquad 
\lim_{t \rightarrow \pm \infty} v_{sph} (t) = 0.$$ 
Frank and K\"onig proved that for each $\epsilon \in (0,v_{cyl}]$ 
there exists a unique positive solution $v_\epsilon$ 
of \eqref{paneitz_ode1} attaining its minimum value of 
$\epsilon$ at $t=0$. Furthermore they prove that any 
positive, global solution of \eqref{paneitz_pde1} must be either 
$v_{sph} (\cdot +T)$ for some $T \in \R$ or $v_\epsilon (\cdot +T)$
for some $T \in \R$ and $\epsilon \in (0,v_{cyl}]$. 

Each $v_\epsilon$ is periodic with period $T_\epsilon$, has local 
minima at exactly $k T_\epsilon$ for each $k \in \mathbf{Z}$, 
and local maxima at exactly $\left ( \frac{2k+1}{2} \right ) T_\epsilon$ 
for each $T_\epsilon$, and no other critical points. For our purposes we 
take the period of the cylindrical solution to be $T_{cyl}$, given in 
\eqref{cyl_period}. Please see the text surrounding \eqref{cyl_period} 
for our reasoning. One can also 
show each $v_\epsilon$ is symmetric about each 
of its critical point and that $T_\epsilon$ is a decreasing function of 
$\epsilon$ with $\lim_{\epsilon \searrow 0} T_\epsilon = \infty$. 
In this context one should think of $v_{sph}$ as the limit 
of $v_\epsilon (\cdot + T_\epsilon /2)$ as $\epsilon \searrow 0$. 

It appears now that the Delaunay metrics occur in a two-parameter 
family, with the necksize $\epsilon$ and translation parameter $T$ as 
parameters. However, it is useful to enlarge this family to include a
$2n$-dimension family of ambient translations, which we first describe 
in Euclidean coordinates and then transform to cylindrical coordinates.
In preparation we transform the Delaunay solution $v_\epsilon$ to 
Euclidean coordinates, obtaining 
$$u_\epsilon(x) = |x|^{\frac{4-n}{2}} v_\epsilon (-\log |x|).$$ 
The first translation is 
\begin{eqnarray*}
\bar u_{\epsilon,a}(x) & = & u_\epsilon(x-a) = |x -a|^{\frac{4-n}{2}}
v_\epsilon( -\log |x-a|) \\
& = & |x|^{\frac{4-n}{2}} \left | \frac{x}{|x|} 
- \frac{a}{|x|} \right |^{\frac{4-n}{2}} v_\epsilon \left (
-\log|x| -\log\left | \frac{x}{|x|} -\frac{a}{|x|} \right | \right ),
\end{eqnarray*}
where $a \in \R^n$. Tranforming this back to cylindrical 
coordinates we then obtain 
\begin{equation} \label{trans_del_soln1}
\bar v_{\epsilon,a} = |\theta - e^t a|^{\frac{4-n}{2}} 
v_\epsilon(t-\log|\theta - e^t a|). \end{equation}
The function $\bar v_{\epsilon,a}$ is of course not a smooth 
global solution, and it has a singular point when 
$$\theta = e^t a \Leftrightarrow t= -\log |a|, \quad \theta 
= \frac{a}{|a|}.$$ 

We obtain the remaining translations using the Kelvin 
transform defined in \eqref{kelvin_trans}. Given $a\in\R^n$ 
we define 
\begin {eqnarray*} 
u_{\epsilon, a} (x) & = & \widehat{\mathbb{K}}_0 (\widehat{\mathbb{K}}_0
(u_\epsilon (\cdot - a)) (x)  =  \widehat{\mathbb{K}}_0 
\left ( |\cdot - a|^{4-n} u_\epsilon \left (
\frac{\cdot - a}{|\cdot - a|^2} \right ) \right ) (x) \\ 
& =& |x|^{4-n} \left | \frac{x}{|x|^2} - a \right |^{4-n}u_\epsilon 
\left ( \frac{\frac{x}{|x|^2} - a}{\left | \frac{x}{|x|^2} - a \right |^2}
\right )  = \left |\frac{x}{|x|} - |x|a\right |^{4-n} u_\epsilon \left ( 
|x|^{-1} \left |\frac{x}{|x|} - |x| a \right |^{-1} \right ) \\ 
& = & \left | \frac{x}{|x|} - |x| a \right |^{\frac{4-n}{2}} |x|^{\frac{n-4}{2}} 
\left | \frac{x}{|x|} - |x| a \right |^{\frac{n-4}{2}} v_\epsilon \left ( 
-\log \left ( |x|^{-1} \left | \frac{x}{|x|} - |x| a \right |^{-1} \right ) 
\right ) \\ 
& = & |x|^{\frac{n-4}{2}} \left | \frac{x}{|x|} - |x| a\right |^{\frac{4-n}{2}}
v_\epsilon \left ( -\log |x| - \log \left | \frac{x}{|x|} - |x| a \right | \right )  ,
\end {eqnarray*}
which in turn gives us
 \begin {equation}\label{trans_del_soln2} 
 v_{\epsilon,a}(t,\theta) = |\theta - e^{-t}a|^{\frac{4-n}{2}} 
 v_\epsilon(t+\log|\theta - e^{-t}a|).\end{equation}
 This function has a singular point when $t=\log|a|$ 
 and $\theta = a/|a|$. 

Using the Taylor expansions 
$$\left | \frac{x}{|x|} - |x| a \right |^{\frac{4-n}{2}} = 
1 + \frac{(n-4)}{2} \langle x,a \rangle + \mathcal{O}(|x|^2) $$
and 
$$\log  \left | \frac{x}{|x|} - a |x| \right | = - \langle a, x \rangle 
+ \mathcal{O} (|x|^2)$$ 
we expand $u_{\epsilon,a}(x)$ as
\begin {eqnarray*} 
u_{\epsilon,a} (x) & = & |x|^{\frac{4-n}{2}} \left ( 1 + \frac{n-4}{2} 
\langle x,a \rangle + \mathcal{O}(|x|^2) \right ) \left (
v_\epsilon(-\log |x|) - \langle a,x \rangle \dot v_\epsilon(-\log |x|) 
+ \mathcal{O}(|x|^2) \right ) \\ 
& = & |x|^{\frac{4-n}{2}} \left ( v_\epsilon(-\log |x|) 
+ \langle x,a\rangle \left ( - \dot v_\epsilon(-\log |x|) + \frac{n-4}{2} 
v_\epsilon(-\log |x|) \right ) + \mathcal{O}(|x|^2) \right ) \\ 
& = & u_\epsilon(x) + |x|^{\frac{4-n}{2}} \langle x,a\rangle 
\left ( -\dot v_\epsilon + \frac{n-4}{2} v_\epsilon \right ) + 
\mathcal{O}(|x|^{\frac{8-n}{2}})  \end {eqnarray*} 
as $|x| \rightarrow 0$. Rewriting this in cylindrical coordinates gives  
\begin{equation} \label{trans_del_soln3} 
v_{\epsilon, a}(t,\theta) = v_\epsilon(t) + e^{-t} \langle \theta, a \rangle 
\left ( - \dot v_\epsilon (t) + \frac{n-4}{2} v_\epsilon(t) \right ) + 
\mathcal{O}(e^{-2t}) 
\end{equation} 
as $t \rightarrow \infty$. 
Unfortunately, the same expansion for $\bar v_{\epsilon,a}$ reveals 
$$\bar v_{\epsilon,a}(t,\theta)  = e^t \left (  \dot v_\epsilon (t) + \frac{n-4}{2}
v_\epsilon (t)  \right )  + \mathcal{O}(1).$$
In fact, 
one expects this behavior, as the motion generating the family 
$\bar v_{\epsilon, a}$ translates the origin in Euclidean coordinates, which 
translates the end in cylindrical coordinates in which $t \rightarrow \infty$. 

It will be convenient for our later computations to observe that 
\begin{eqnarray} \label{trans_nodal_domains1} 
&\langle a,\theta\rangle  > 0 \Rightarrow v_{\epsilon,a} (t,\theta) > v_\epsilon(t) , \qquad 
\langle a,\theta \rangle < 0 \Rightarrow v_{\epsilon, a} (t,\theta) < v_\epsilon(t) & 
\\ \nonumber 
& \langle a,\theta \rangle > 0 \Rightarrow \bar v_{\epsilon,a} (t,\theta) > 
v_\epsilon (t) , \qquad \langle a,\theta \rangle < 0 \Rightarrow \bar v_{\epsilon,a} 
(t,\theta) < v_\epsilon (t) .& 
\end{eqnarray} 

One can find a first integral for the ODE \eqref{paneitz_ode1}. Indeed, 
differentiating once shows 
\begin{equation} \label {delaunay_energy}
\mathcal{H}_\epsilon = -\dot v_\epsilon \dddot v_\epsilon
+ \frac{1}{2} \ddot v_\epsilon^2 + \left ( \frac{n(n-4) + 8}{4} \right ) 
\dot v_\epsilon^2 - \frac{n^2(n-4)^2}{32} v_\epsilon^2 + \frac{(n-4)^2
(n^2-4)}{32} v_\epsilon^{\frac{2n}{n-4}} \end{equation} 
is constant function for each Delaunay solution $v_\epsilon$. 
Evaluating this energy on the cylindrical 
solution gives 
\begin{equation} \label{cyl_energy} 
\mathcal{H}_{cyl}  = -\frac{(n-4)(n^2-4)}{8} \left ( \frac{n(n-4)}{n^2-4} 
\right )^{\frac{n}{4}} < 0 \end{equation}
and evaluating on the spherical solution gives 
\begin{eqnarray} \label{sph_energy} 
\mathcal{H}_{sph} & = & \frac{(n-4)^2}{4} (\cosh t)^{-n} \left ( 
\frac{n(n-2)}{4} \sinh^4 t+ \left ( \frac{4-n}{2} + 2-n\right ) \cosh^2 t \sinh^2 t 
- \frac{1}{2} \cosh^4 t \right . \\ \nonumber 
&& \qquad \qquad + \frac{(n-2)}{4}\cosh^2 t\sinh^2 t - \frac{(n-2)^2}{8} 
\sinh^4 t - \left ( \frac{n(n-4)+8}{4} \right ) \cosh^2 t\sinh^2 t \\ \nonumber 
&& \qquad \qquad \left. + \frac{n^2}{8} 
\cosh^4 t - \frac{(n^2-4)}{8} \right ) \\ \nonumber 
& = & \frac{(n-4)^2}{4}(\cosh t)^{-n} \left ( -\frac{(n^2-4)}{8} + \frac{(n^2-4)}{8}
\cosh^4 t + \frac{(4-n^2)}{4} \cosh^2 t \sinh^2 t + \frac{(n^2-4)}{8} 
\sinh^4 t \right ) \\ \nonumber 
& = & 0. 
\end{eqnarray}

It is a consequence of Proposition 6 of \cite{FK} that $\mathcal{H}_\epsilon$
is a strictly decreasing function of $\epsilon$. 

\subsection{Linearization of the PDE \eqref{paneitz_pde2}} 

Let $u$ solve \eqref{paneitz_pde1} and consider the 
slightly perturbed function $u+ \epsilon w$. We have 
\begin {equation} \label{linear_rhs1}
\frac{n(n-4)(n^2-4)}{16} (u+\epsilon w)^{\frac{n+4}{n-4}} 
= \frac{n(n-4)(n^2-4)}{16} u^{\frac{n+4}{n-4}} + \epsilon 
\frac{n(n+4)(n^2-4)}{16} w u^{\frac{8}{n-4}} + \mathcal{Q}(\epsilon w), 
\end {equation} 
where 
\begin {equation} \label {linear_rhs2}
\mathcal{Q}(\epsilon w) = \frac{n(n-4)(n^2-4)}{16} \left [ 
(u+\epsilon w)^{\frac{n+4}{n-4}} - u^{\frac{n+4}{n-4}} 
- \epsilon \frac{n+4}{n-4} w u^{\frac{8}{n-4}} \right ] = \mathcal{O}(\epsilon^2
\| w\|^2).\end {equation} 
Combining \eqref{paneitz_pde1} with \eqref{linear_rhs1} and 
\eqref{linear_rhs2} we have 
$$\Delta^2 u + \epsilon \Delta^2 w = \frac{n(n-4)(n^2-4)}{16} 
u^{\frac{n+4}{n-4}} + \epsilon \frac{n(n+4)(n^2-4)}{16} w 
u^{\frac{8}{n-4}}
+ \mathcal{Q}(\epsilon w),$$
where $\mathcal{Q}(\epsilon w)$ is of 
order $\epsilon^2 \| w \|^2$. 
Selecting the terms of order $\epsilon$ from the last equation 
we obtain the linearization 
\begin {equation}\label{linearization1}
\Delta^2 w = \frac{n(n+4)(n^2-4)}{16} u^{\frac{8}{n-4}} w.
\end {equation}
We refer to a solution $w$ of \eqref{linearization1} as a 
{\bf Jacobi field} associated to the solution $u$. 

Finally, it is useful to rewrite this last linearized equation in 
cylindrical coordinates, obtaining 
$$\frac{\partial^4 w}{\partial t^4} - \left ( \frac{n(n-4)
+8}{2} \right )\frac{\partial^2 w}{\partial t^2} + \frac{n^2(n-4)^2}
{16} w + 2 \Delta_\theta \frac{\partial^2 w}{\partial t^2} 
- \frac{n(n-4)}{2} \Delta_\theta w + \Delta_\theta^2 w 
= \frac{n(n+4)(n^2-4)}{16} w v^{\frac{8}{n-4}},$$
which we can rearrange to read 
\begin {eqnarray} \label{linearization2} 
0 & = & \frac{\partial^4 w}{\partial t^4} - \left ( \frac{n(n-4)
+ 8}{2} \right ) \frac{\partial^2 w}{\partial t^2} + 
\left ( \frac{n^2(n-4)^2}{16} - \frac{n(n+4)(n^2-4)}{16} 
v^{\frac{8}{n-4}} \right ) w \\ \nonumber 
&& + 2 \Delta_\theta \frac{\partial^2 w}{\partial t^2} 
- \frac{n(n-4)}{2} \Delta_\theta w + \Delta_\theta^2 w. 
\end {eqnarray}
Again, we refer to a solution $w$ of \eqref{linearization2} as 
a Jacobi field associated to the solution $v$. 
  
\subsection {Integral identities}

The following is essentially a special case of Proposition 4.2 in \cite{DMO}, 
and also a special case of Proposition A.2 of \cite{JX}. 
\begin {prop} \label{poho_prop}
Let $u \in \mathcal{C}^\infty((0,\infty) \times \mathbf{S}^{n-1})$ 
solve 
\eqref{paneitz_pde2} and let $0<T_1< T_2$. Then 
\begin {eqnarray} \label{poho1} 
& \int_{\{ T_1\} \times \mathbf{S}^{n-1}} \frac{\partial v}
{\partial t} \frac{\partial^3 v}{\partial t^3} - \frac{1}{2} 
\left ( \frac{\partial^2 v}{\partial t^2} \right )^2 - \left ( 
\frac{n(n-4) + 8}{4} \right ) \left (\frac{\partial v}{\partial t}
\right )^2 + \frac{n^2(n-4)^2}{32} v^2 & \\ \nonumber 
& \qquad - \frac{(n-4)^2(n^2-4)}{32} v^{\frac{2n}{n-4}} + 
\frac{1}{2} (\Delta_\theta v)^2 + \frac{n(n-4)}{4} |\nabla_\theta 
v|^2 - \left | \nabla_\theta \frac{\partial v}{\partial t} 
\right |^2 d\sigma \\ \nonumber
= & \int_{\{ T_2\} \times \mathbf{S}^{n-1}} \frac{\partial v}
{\partial t} \frac{\partial^3 v}{\partial t^3} - \frac{1}{2} 
\left ( \frac{\partial^2 v}{\partial t^2} \right )^2 - \left ( 
\frac{n(n-4) + 8}{4} \right ) \left (\frac{\partial v}{\partial t}
\right )^2 + \frac{n^2(n-4)^2}{32} v^2 & \\ \nonumber 
& \qquad - \frac{(n-4)^2(n^2-4)}{32} v^{\frac{2n}{n-4}} + 
\frac{1}{2} (\Delta_\theta v)^2 + \frac{n(n-4)}{4} |\nabla_\theta 
v|^2 - \left | \nabla_\theta \frac{\partial v}{\partial t} 
\right |^2 d\sigma \end {eqnarray} 
\end {prop} 

\begin {proof} Multiplying \eqref{paneitz_pde2} by 
$\frac{\partial v}{\partial t}$ and integrating over the sphere 
$\{ t \} \times \mathbf{S}^{n-1}$ we obtain 
\begin {eqnarray*} 
0 & = & \int_{\{ t\} \times \mathbf{S}^{n-1}}
\frac{\partial v}{\partial t} \frac{\partial^4 v}{\partial 
t^4} - \frac{n(n-4) + 8}{2} \frac{\partial^2 v}{\partial t^2} 
\frac{\partial v}{\partial t} + \frac{n^2(n-4)^2}{16} v \frac
{\partial v}{\partial t} \\ 
&& \qquad + \frac{\partial v}{\partial t} \Delta^2_\theta v- 
\frac{n(n-4)}{2} \frac{\partial v}{\partial t} \Delta_\theta v
+ 2 \frac{\partial v}{\partial t} \Delta_\theta \frac{\partial^2 v}
{\partial t^2} - \frac{n(n-4)(n^2-4)}{16} v^{\frac{n+4}{n-4}}
\frac{\partial v}{\partial t} d\sigma \\ 
& = & \int_{\{ t \} \times \mathbf{S}^{n-1}} \frac{\partial}{\partial t} 
\left [ \frac{\partial v}{\partial t}\frac{\partial^3 v}{\partial t^3} 
- \frac{1}{2} \left ( \frac{\partial^2 v}{\partial t^2} \right )^2
- \left ( \frac{n(n-4) + 8}{4} \right ) \left (\frac{\partial v}
{\partial t} \right )^2 + \frac{n^2(n-4)^2}{32} v^2 \right . \\ 
& & \qquad \qquad - \frac{(n-4)^2(n^2-4)}{32} v^{\frac{2n}
{n-4}} - \frac{1}{2} (\Delta_\theta v)^2 + \frac{n(n-4)}{4} 
|\nabla_\theta v|^2 \\
&& \qquad \qquad - \left .  \frac{n(n-4)}{2} \operatorname{div}_\theta
\left ( \frac{\partial v}{\partial t} \nabla_\theta v \right ) 
- \left | \nabla_\theta \frac{\partial v}{\partial t} \right |^2 
+ 2 \operatorname{div}_\theta \left ( \frac{\partial v}{\partial t} 
\nabla_\theta \frac{\partial^2 v}{\partial t^2} \right ) \right ] d\sigma  \\ 
& = & \frac{d}{dt} \int_{\{ t \} \times \mathbf{S}^{n-1}}
\frac{\partial v}{\partial t}\frac{\partial^3 v}{\partial t^3} 
- \frac{1}{2} \left ( \frac{\partial^2 v}{\partial t^2} \right )^2
- \left ( \frac{n(n-4) + 8}{4} \right ) \left (\frac{\partial v}
{\partial t} \right )^2 + \frac{n^2(n-4)^2}{32} v^2\\ 
& & \qquad \qquad - \frac{(n-4)^2(n^2-4)}{32} v^{\frac{2n}
{n-4}} - \frac{1}{2} (\Delta_\theta v)^2 + \frac{n(n-4)}{4} 
|\nabla_\theta v|^2 
- \left | \nabla_\theta \frac{\partial v}{\partial t} \right |^2 
 d\sigma. 
\end{eqnarray*}
\end {proof} 

Following this integral identity, we define 
\begin{eqnarray} \label{radial_poho}
\mathcal{H}_{rad}(v) & = & \int_{ \{ t\} \times \Ss^{n-1}} 
-\frac{\partial v}{\partial t} \frac{\partial^3 v}{\partial t^3} + 
\frac{1}{2} \left ( \frac{\partial^2 v}{\partial t^2} \right )^2 + 
\left ( \frac{n(n-4)+8}{4} \right )\left ( \frac{\partial v}{\partial t} \right )^2 
\\ \nonumber  
&& \quad - \frac{n^2(n-4)^2}{32} v^2 
+ \frac{(n-4)^2(n^2-4)}{32} v^{\frac{2n}{n-4}} - \frac{1}{2} 
(\Delta_\theta v)^2 - \frac{n(n-4)}{4} |\nabla_\theta v|^2 + \left | 
\nabla_\theta \frac{\partial v}{\partial t} \right |^2 d\sigma ,
\end{eqnarray} 
which does not depend on $t$ by \eqref{poho1}. 

\begin {cor} 
If $v_\epsilon$ is the Delaunay solution of \eqref{paneitz_ode1}
described above then 
\begin{equation} \label{poho2} 
\mathcal{H}_{rad} (v_\epsilon) 
= n \omega_n \left ( \frac{1}{2} \ddot v_\epsilon^2 (0)+ \frac{(n-4)^2}{32} 
\epsilon^2 [(n^2-4) \epsilon^{\frac{8}{n-4}} - n^2] \right ) 
=n \omega_n \mathcal{H}_\epsilon  < 0 .
\end{equation}
where $\omega_n$ is the volume of a unit ball in $\R^n$. 
\end {cor} 

\subsection{A priori estimates} 

Jin and Xiong \cite{JX} prove the following. 
\begin{thm} Let $u\in \mathcal{C}^\infty (\B_1(0) \backslash 
\{ 0 \})$ be a positive solution of \eqref{paneitz_pde1} such 
that $-\Delta u > 0$. Then either $u$ admits a continuous extension 
to $\B_1(0)$ or there exist positive constants $C_1< C_2$ 
such that 
\begin{equation} \label{apriori_bounds2}
C_1 |x|^{\frac{4-n}{2}} \leq u(x) \leq C_2 |x|^{\frac{4-n}{2}}. 
\end{equation}
\end{thm} 

\begin {rmk} The example of the Delaunay solutions demonstrate that the 
constant $C_1$ in \eqref{apriori_bounds2} must depends on the choice of the 
solution $u$. However, the constant $C_2$ in the upper bound is universal. 
\end {rmk}

\begin{lemma} Let $v \in \mathcal{C}^\infty ((0, \infty) \times \Ss^{n-1})$ 
be a positive solution of \eqref{paneitz_pde2} which also 
satisfies \eqref{pos_lap_cyl_coords}. Then either $\limsup_{t 
\rightarrow \infty} v(t,\theta) = 0$ or there exist positive 
constants $C_1 < C_2$ such that 
$$C_1 < v(t,\theta) < C_2 . $$ 
Moreover, in this case $\mathcal{H}_{rad} (v) \neq 0$. 
\end{lemma} 

\begin{proof} The bounds \eqref{apriori_bounds2} imply 
$C_1 < v(t,\theta)< C_2$. Choose a sequence $\tau_i \nearrow 
\infty$ and define $v_i (t,\theta) = v(t+\tau_i, \theta)$. This 
sequence is uniformly bounded, so we may extract a 
subsequence, still denoted by $\{ v_i\}$, which converges 
uniformly on compact subsets of $\R \times \Ss^{n-1}$
to a solution $\bar v$ 
of \eqref{paneitz_pde2}. However, the limit $\bar v$ must be a 
Delaunay solution $v_\epsilon$, and so (using  Proposition \ref{poho_prop}) 
we have 
\begin{eqnarray*} 
\mathcal{H}_{\textrm{rad}} (v) & = & \int_{ \{ 1 \} \times \Ss^{n-1}} 
- \frac{\partial v}{\partial t} 
\frac{\partial^3 v}{\partial t^3}+ \frac{1}{2} \left ( \frac{\partial^2 v}
{\partial t^2} \right )^2 + \left ( \frac{n(n-4) + 8}{4} \right ) \left ( \frac 
{\partial v}{\partial t} \right )^2 \\ 
&& \quad + \frac{n^2(n-4)^2}{32} v^2 + \frac{(n-4)^2(n^2-4)}{32} 
v^{\frac{2n}{n-4}} - \frac{1}{2} (\Delta_\theta v)^2- \frac{n(n-4)}{4}
\left | \nabla_\theta \frac{\partial v}{\partial t} \right |^2 d\sigma \\ 
 &= & \int_{ \{ \tau_i +1 \} \times \Ss^{n-1}} 
- \frac{\partial v}{\partial t} 
\frac{\partial^3 v}{\partial t^3}+ \frac{1}{2} \left ( \frac{\partial^2 v}
{\partial t^2} \right )^2 + \left ( \frac{n(n-4) + 8}{4} \right ) \left ( \frac 
{\partial v}{\partial t} \right )^2 \\ 
&& \quad + \frac{n^2(n-4)^2}{32} v^2 + \frac{(n-4)^2(n^2-4)}{32} 
v^{\frac{2n}{n-4}} - \frac{1}{2} (\Delta_\theta v)^2- \frac{n(n-4)}{4}
\left | \nabla_\theta \frac{\partial v}{\partial t} \right |^2 d\sigma \\ 
& = & \mathcal{H}_{\textrm{rad}} (v_i) \rightarrow 
\mathcal{H}_{\textrm{rad}} (v_\epsilon) = n \omega_n \mathcal{H}_\epsilon 
\neq 0.
\end{eqnarray*} 
\end{proof}

\section{Linear analysis}
\label{lin_anal_sec}

In this section we study the mapping properties of the linear 
operator \eqref{linearization2}, concentrating on the linearization 
about a Delaunay metric. 

\subsection{Definitions}

Linearizing \eqref{paneitz_pde2} about a Delaunay solution $v_\epsilon$ 
we obtain the operator 
\begin {eqnarray} \label{del_linearization1}
L_\epsilon & = & \frac{\partial^4}{\partial t^4} + \Delta_\theta^2 + 2 \Delta_\theta 
\frac{\partial^2}{\partial t^2} - \frac{n(n-4)}{2} \Delta_\theta \\ \nonumber 
&& - \left ( \frac{n(n-4) + 8} 
{2} \right ) \frac{\partial^2 }{\partial t^2} + \left ( \frac{n^2 (n-4)^2}{16} - 
\frac{n(n+4)(n^2-4)}{16} v_\epsilon^{\frac{8}{n-4}} \right ) .
\end{eqnarray}
We refer to solutions of the equation $L_\epsilon (w)$ as {\bf Jacobi 
fields} associated to the Delaunay solution $v_\epsilon$. More generally, 
if $v$ satisfies \eqref{paneitz_pde2} and $w$ satisfies \eqref{linearization2}
then we call $w$ a Jacobi field associated to the solution $v$. 
We are interested in the mapping properties of $L_\epsilon$, for 
instance as the map 
$$L_\epsilon : W^{4,2}((0,\infty) \times \Ss^{n-1}) \rightarrow L^2
((0, \infty) \times \Ss^{n-1}).$$
It turns out the operator written above does not have closed 
range (see pg. 21, pg. 216, and Theorem 5.40 of \cite{Mel}), so we 
will to define certain weighted function spaces to 
accomodate $L_\epsilon$.

\begin {defn}
Let $\gamma \in \R$ and when 
$u \in L^2((0,\infty) \times \Ss^{n-1})_{loc}$ define 
\begin {equation} \label{cyl_weight_norm}
\| u\|^2_{L^2_\gamma} = \int_0^\infty \int_{\Ss^{n-1}} 
e^{-2\gamma t} | u(t, \theta)|^2 d\theta dt .
\end {equation}
The space $L^2_\gamma((0,\infty) \times \Ss^{n-1})$ is the 
space of functions with finite norm, as defined above. One can 
similarly define the Sobolev spaces $W^{k,2}_\gamma$ of functions 
with $k$ weak derivatives in $L^2$ having finite weighted norms.
\end {defn} 
One can also define weighted H\"older spaces. 
\begin {defn} Let $\gamma\in \R$ and $\alpha \in (0,1)$. For 
$u \in \mathcal{C}^{0,\alpha}_{loc}((0,\infty) \times \Ss^{n-1})$ define 
$$\| u \|_{\mathcal{C}^{0,\alpha}_\gamma} =  
\sup_{t_0 > 1} \sup \left \{ \frac{ e^{-\gamma t_1} u(t_1,\theta_1) - 
e^{-\gamma t_2} u(t_2, \theta_2)}{d((t_1, \theta_1), (t_2, \theta_2))^\alpha} :
(t_1, \theta_1), (t_2, \theta_2) \in (t_0-1, t_0+1) \times \Ss^{n-1} \right \}.$$
One can similarly define weighted H\"older spaces with more 
derivatives. 
\end {defn}
Heuristically, a function in a weighted function space with weight $\gamma$ 
is bounded from above by a multiple of $e^{\gamma t}$ as $t \rightarrow \infty$.
Observe then that when $\gamma_1 < \gamma_2$ we have the inclusions 
$$W^{k,2}_{\gamma_1}((0,\infty) \times \Ss^{n-1}) \subset 
W^{k,2}_{\gamma_2}((0,\infty) \times \Ss^{n-1}), \qquad 
\mathcal{C}^{k,\alpha}_{\gamma_1}((0,\infty) \times \Ss^{n-1}) \subset 
\mathcal{C}^{k,\alpha}_{\gamma_2}((0,\infty) \times \Ss^{n-1}) .$$

The fact that $L_\epsilon$ the leading order terms of $L_\epsilon$ is 
$\frac{\partial ^4}{\partial t^4} + \Delta_\theta^2$ tells us the following. 
\begin {lemma} 
For any $\gamma \in \R$ the operators
$$L_\epsilon : W^{k+4, 2}_\gamma((0, \infty) \times \Ss^{n-1}) \rightarrow W^{k,2}_\gamma
((0, \infty) \times \Ss^{n-1}) $$
and 
$$L_\epsilon : \mathcal{C}^{k+4, \alpha}_\gamma((0, \infty) \times \Ss^{n-1}) \rightarrow 
\mathcal{C}^{k,\alpha}_\gamma ((0, \infty) \times \Ss^{n-1})$$
are bounded, linear, and elliptic. \end {lemma} 

We complete our understanding of $L_\epsilon$ by identifying the 
weights for which 
$$L_\epsilon : W^{k+4, 2}_\gamma((0, \infty) \times \Ss^{n-1}) \rightarrow W^{k,2}_\gamma
((0, \infty) \times \Ss^{n-1}) $$
is Fredholm, injective, and/or surjective. Indeed, this is a nontrivial 
task. 

We make our analysis easier by decomposing $w$ in 
spherical harmonics, writing 
$$w(t,\theta) = \sum_{j \in \mathbf{Z}} w_j (t) \phi_j,$$
where $\phi_j$ is a normalized eigenfunction of $\Delta_\theta$ 
on $\Ss^{n-1}$, {\it i.e.} 
\begin{equation} \label{spherical_eigenfunc} 
\Delta_\theta \phi_j = - \lambda_j \phi_j, \qquad \int_{\Ss{n-1}} 
\phi_j \phi_k d\theta = \delta_{jk}.
\end{equation} 
The eigenvalues $\lambda_j$ of the $(n-1)$-dimensional
sphere have the form $\lambda_j = k(n-2+k)$ for 
some $k=0,1,2,3,\dots$.

Under this decomposition the Fourier coefficient $w_j$ solves the 
ODE 
\begin{eqnarray}  \label{paneitz_ode3}
0 & = & L_{\epsilon, j} w_j \\ \nonumber 
& = & \ddddot w_j - \left ( \frac{n(n-4) + 8 + 4\lambda_j}{2} \right ) \ddot w_j 
+ \left ( \frac{n^2(n-4)^2}{16} - \frac{n(n+4)(n^2-4)}{16} v_\epsilon^{\frac{8}{n-4}}
+ \frac{n(n-4)}{2} \lambda_j+ \lambda_j^2 \right ) w_j. 
\end{eqnarray} 
It immediately follows that 
\begin{equation} \label{spec_decomp}
\operatorname{spec}(L_\epsilon) = \bigcup_{j=0}^\infty \operatorname{spec}
(L_{\epsilon,j}).\end{equation} 

\subsection{Low Fourier modes}

We can explicitly identify some of the ODE solutions when $|j|$ is small. 
For instance, we have $\lambda_0 = 0$ and so \eqref{paneitz_ode3} becomes
\begin {equation} \label{paneitz_ode4}
\ddddot w_0 - \left ( \frac{n(n-4) +8}{2} \right ) \ddot w_0 + 
\left ( \frac{n^2(n-4)^2}{16} - \frac{n(n+4)(n^2-4)}{16} v_\epsilon^{\frac{8}
{n-4}} \right )w_0 =0,\end{equation}
which is the derivative of \eqref{paneitz_ode1}, Thus 
\begin{equation} \label{low_paneitz_mode1}
w_0^+ (t) = \dot v_\epsilon(t) , \qquad w_0^- = \frac{d}{d\epsilon} 
v_\epsilon (t) 
\end{equation}
both solve \eqref{paneitz_ode4}.  
\begin {lemma} \label{zero_ind_roots}
The function $w_0^+$ is periodic with 
period $T_\epsilon$ while the function $w_0^-$ grows 
linearly. \end {lemma}

\begin {proof} 
Differentiating the equation $v_\epsilon(t) = v_\epsilon(t+T_\epsilon)$ 
with respect to $t$ gives $w_0^+(t+T_\epsilon) = w_0(t)$. Differentiating 
$v_\epsilon(t+T_\epsilon) = v_\epsilon(t)$ with respect to $\epsilon$ 
gives $$w_0^-(t+T_\epsilon) \frac{d T_\epsilon}{d\epsilon} = w_0^-
(t),$$ and so $w_0^-$ grows linearly. \end {proof} 

We can also explicitly identify the Fourier modes when 
$j=1,2,\dots, n$ and $\lambda_j = n-1$. To do this we first 
let $e_j$ be the standard basis of $\R^n$ and observe 
$\phi_j = \langle e_j, \theta \rangle$ is the eigenfunction 
associated to $\lambda_j = n-1$. Substituting $a= \tau e_j$ 
into \eqref{trans_del_soln2} and \eqref{trans_del_soln3} we 
find 
\begin{eqnarray*} 
v_{\epsilon, \tau e_j}  & = & |\theta - e^{-t} \tau e_j |^{\frac{4-n}{2}} 
v_\epsilon (t +\log |\theta - e^{-t} \tau \theta|) \\ 
& = & v_\epsilon(t) + \tau e^{-t} \langle \theta, e_j \rangle \left ( 
-\dot v_\epsilon (t) + \frac{n-4}{2} v_\epsilon (t) \right ) + \mathcal{O} 
(e^{-2t}) \\ 
& = & v_\epsilon (t) + \tau e^{-t} \phi_j \left ( -\dot v_\epsilon (t) 
+ \frac{n-4}{2} v_\epsilon (t) \right ) + \mathcal{O} (e^{-2t}). 
\end{eqnarray*} 
Differentiating with respect to $\tau$ we obtain 
$$w_j^- (t) \phi_j(\theta) = \left. \frac{d}{d\tau} \right |_{\tau = 0} 
v_{\epsilon, \tau e_j} = e^{-t} \left ( -\dot v_\epsilon (t) + \frac{n-4}{2} 
v_\epsilon (t) \right ) \phi_j(\theta) + \mathcal{O}(e^{-2t}) ,$$
or 
\begin {equation} \label{low_paneitz_mode2} 
w_j^-  = e^{-t} \left ( -w_0^+ + \frac{n-4}{2} v_\epsilon \right ) 
+ \mathcal{O}(e^{-2t}) = e^{-t}  \left (- \dot v_\epsilon + \frac{n-4}{2} v_\epsilon 
\right ) + \mathcal{O} (e^{-2t} ). \end {equation} 

However, $v_{\epsilon, \tau e_j}$ satisfies \eqref{paneitz_pde2} 
for each $\tau$. Differentiating this relation and using $-\Delta_\theta 
\phi_j = (n-1)\phi_j = \lambda_j \phi_j$ we find 
\begin{eqnarray*} 
\frac{n(n+4)(n^2-4)}{16} v_\epsilon^{\frac{8}{n-4}} w_j^- \phi_j & = & 
\left ( \frac{\partial^4}{\partial t^4} - \left ( \frac{n(n-4)+8}{2} \right ) \frac{\partial^2}
{\partial t^2} + \frac{n^2(n-4)^2}{16} \right ) w_j^- \phi_j \\ 
&& + \left ( \Delta_\theta^2 + 2 \Delta_\theta \frac{\partial^2}{\partial t^2} - 
\frac{n(n-4)}{2} \Delta_\theta \right ) w_j^- \phi_j \\ 
& = & \left ( \ddddot w_j^- - \left ( \frac{n(n-4)+8+4\lambda_j}{2} \right ) \ddot w_j^- 
+ \left ( \frac{n^2 (n-4)^2}{16} + \lambda_j^2 + \frac{n(n-4)}{2} \lambda_j \right ) 
w_j^-  \right ) \phi_j, 
\end{eqnarray*} 
which we can rearrange to give $L_{\epsilon, j} (w_j^-) = 0$. 

A similar calculation, starting from 
$$w_j^+ (t) \phi_j(\theta) = \left. \frac{d}{d\tau} \right |_{\tau = 0}
\bar v_{\epsilon, \tau e_j}, $$
gives the expansion 
\begin {equation} \label{low_paneitz_mode3}  
w_j^+  = e^t \left ( w_0^+ + \frac{n-4}{2} v_\epsilon  \right ) 
+ \mathcal{O}(1) = e^t \left ( \dot v_\epsilon + \frac{n-4}{2} 
v_\epsilon \right ) + \mathcal{O}(1)  .\end {equation} 

Observe that \eqref{trans_nodal_domains1} implies 
\begin{equation} \label{trans_nodal_domains2} 
w_j^+ > 0, \qquad w_j^- > 0 . 
\end{equation} 

It is not surprising that $w_j^\pm(t)$ all agree for $j=1,2,\dots, n$, 
as each translation is geometrically the same. It is also not 
surprising that $w_j^+$ grows exponentially while $w_j^-$ 
decays exponentially. The ambient motion generating $w_j^+$ 
translates the origin in Euclidean coordinates, which in cylindrical 
coordinates moves the end corresponding to $t \rightarrow \infty$, 
whereas the ambient motion generating $w_j^-$ 
moves the end corresponding to $t\rightarrow -\infty$.  

\subsection{Indicial roots} 
\label{indicial_sect}

We further analyze the ODE \eqref{paneitz_ode3} for general values 
of $j$, that is 
\begin{eqnarray*}
0 & = & L_{\epsilon, j} w \\
& = & \ddddot w - \left ( \frac{n(n-4) + 8 + 4\lambda_j}{2} \right ) \ddot w + 
\left ( \frac{n^2(n-4)^2}{16} - \frac{n(n+4)(n^2-4)}{16} v_\epsilon^{\frac{8}{n-4}}
+\lambda_j^2 + \frac{n(n-4)}{2} \lambda_j \right ) w .\end{eqnarray*}
The coefficients of this ODE are all periodic with period $T_\epsilon$, so 
there is a (constant) $4 \times 4$ matrix $A_{\epsilon,j}$ such that 
\begin{equation}\label{indicial1}
\left ( \begin {array}{c} w(t+T_\epsilon)\\ \dot w(t+T_\epsilon) \\ 
\ddot w(t+T_\epsilon) \\ \dddot w(t+T_\epsilon) \end {array} \right ) = 
A_{\epsilon,j} \left ( \begin {array}{c} w(t) \\ \dot w(t) \\ \ddot w(t) \\ 
\dddot w(t) \end {array} \right ).\end{equation} 
Now let $w_1, w_2,w_3,w_4$ be four solutions of \eqref{paneitz_ode3}, and 
let 
\begin{equation}\label{wronskian1} 
W(t) = \det \left ( \begin {array}{cccc} w_1 & w_2 & w_3 & w_4 \\ 
\dot w_1 & \dot w_2 & \dot w_3 & \dot w_4 \\ \ddot w_1 & \ddot w_2
& \ddot w_3 & \ddot w_4 \\ \dddot w_1 & \dddot w_2 & \dddot w_3 
& \dddot w_4 \end {array} \right ) (t) \end{equation} 
be the associated Wronskian determinant. By Abel's identity, 
$\frac{dW}{dt} = 0$ and so $W$ is constant. Combining \eqref{indicial1} 
and \eqref{wronskian1} we see that $\det A_{\epsilon,j} = 1$. Moreover, 
the matrix $A_{\epsilon,j}$ has real coefficients, so its eigenvalues occur 
in conjugate pairs. Suppressing the dependence on $\epsilon$ 
for the moment, we denote these eigenvalues as 
\begin{equation}\label{indicial2}
\mu_j^\pm = e^{\pm i \xi_{\epsilon,j}}, \qquad 
\tilde \mu_j^\pm = e^{\pm i \tilde \xi_{\epsilon,j}} .
\end{equation} 
Here we define the {\bf indicial roots} of the operator 
$L_{\epsilon,j}$ as all the real numbers $\gamma$ such that 
$\gamma = \Im (\xi)$ where $\mu=e^{i\xi}$ is an eigenvalue of $A_{\epsilon,j}$. 
For convenience later on, we collect these numbers as 
\begin{equation} \label{indicial_defn} 
\Gamma_{\epsilon,j} = \{ \gamma  \in \R :
\gamma = \Im(\xi) \textrm{ and } \mu =e^{i\xi}
\textrm{ is an eigenvalue of }A_{\epsilon,j} \}, \qquad 
\Gamma_\epsilon = \bigcup_{j=0}^\infty \Gamma_{\epsilon,j}. 
\end{equation}
Observe that $\Gamma_{\epsilon,j}$ has at most four elements, so 
in particular $\Gamma_\epsilon$ is countable. Moreover, 
by construction $\Gamma_{\epsilon,j}$ is even for each $j$, 
{\it i.e.} $\gamma \in \Gamma_{\epsilon, j}$ if and only if 
$-\gamma \in \Gamma_{\epsilon,j}$. 

The fact that the eigenvalues of $A_{\epsilon,j}$ occur in 
conjugate pairs implies we can always write the eigenvalues 
as $\mu^\pm = e^{\pm i \xi}$ and $\tilde \mu^\pm = 
e^{\pm i \tilde \xi}$. Writing $\xi = \eta + i \nu$ and $\tilde 
\xi = \tilde \eta + i \tilde \nu$ we see that $|\mu^\pm| = e^{\pm 
\nu}$ and $|\tilde \mu^\pm| = e^{\pm \tilde \nu}$. In other 
words, the indicial roots precisely 
determine the rates of exponential growth of the solutions of 
\eqref{paneitz_ode3}. 

\begin{lemma} 
For each $\epsilon \in (0,v_{cyl}]$ we have $0 \in 
\Gamma_{\epsilon,0}$ (with multiplicity $2$) and $\{ 
-1,1 \} \subset \Gamma_{\epsilon,j}$ for $j=1,2,\dots, n$. 
\end{lemma} 
\begin{proof} Lemma \ref{zero_ind_roots} implies 
$0 \in \Gamma_{\epsilon,0}$ while \eqref{low_paneitz_mode2}
and \eqref{low_paneitz_mode3} together imply $\{ \pm 1\} 
\subset \Gamma_{\epsilon, j}$ for $1\leq j \leq n$. 
\end{proof} 

Now we explicitly compute the Jacobi fields and the 
indicial roots for the special case of the 
cylindrical metric. We let $\epsilon_n$ denote the Delaunay 
parameter of the cylindrical solution, that is 
\begin{equation} \label{cyl_eps}
\epsilon_n = \left ( \frac{n(n-4)}{n^2-4} \right ) ^{\frac{n-4}{8}}, 
\end{equation} 
so that 
\begin{equation} \label{cyl_fourier_odes}
L_{\epsilon_n, j} = \frac{d^4}{dt^4} - \left ( \frac{n(n-4) + 8  
+ 4\lambda_j}{2} \right ) \frac{d^2}{dt^2} + \left ( -\frac{n^2(n-4)}{2}  
+ \frac{n(n-4)}{2} \lambda_j + \lambda_j^2 \right ) .
\end{equation} 
Fortunately we can explicity solve the ODEs $L_{\epsilon_n, j} w = 0$. 
Substituting 
$$w(t) = c_+ e^{\mu_j t} + c_- e^{-\mu_j t} + \widetilde c_+ 
e^{\widetilde \mu_j t} + \widetilde c_- e^{-\widetilde \mu t} $$ 
and using $\lambda_j = k(n-2+k)$ for some nonnegative integer $k$ 
we find 
\begin{eqnarray} \label{cyl_fourier_exp1}
\mu_j^2 & = &  \frac{1}{2} \left ( \frac{n(n-4)+8+4\lambda_j}{2} + 
\sqrt{\frac{n^4}{4} - 16(n-1-\lambda_j)} \right ) \\ \nonumber 
& = & \frac{1}{2} \left ( \frac{(n+2(k-1))^2+ 4}{2} + \sqrt{ \frac{n^4}{4} 
+ 16(k-1)(n+k-1)} \right )
\end{eqnarray}
and 
\begin{eqnarray} \label{cyl_fourier_exp2}  
\widetilde \mu_j^2 & = & \frac{1}{2} \left ( \frac{n(n-4)+8+4\lambda_j}{2} -
\sqrt{\frac{n^4}{4} - 16(n-1-\lambda_j)} \right )  \\ \nonumber 
& = &  \frac{1}{2} \left ( \frac{(n+2(k-1))^2+ 4}{2} - \sqrt{ \frac{n^4}{4} 
+ 16(k-1)(n+k-1)} \right ) , 
\end{eqnarray}
which then gives all solutions after taking square roots. 
We remark on some properties of $\mu_j$ and $\widetilde 
\mu_j$. First observe that, because $n > 4$,  
$$\frac{n^4}{4} - 16n +16+16\lambda_j \geq n^4-16n +16 > 16\left ( 
\frac{n^2}{4} -n+1 \right ) >0,$$
so that, in particular, $\mu_j^2$ and $\widetilde \mu_j^2$ are real 
numbers. Next we observe that $\mu_j^2 > 0$ for each integer $j$. 
After taking a positive 
and a negative square root gives us one Jacobi field which grows exponentially 
and one exponentially decaying Jacobi field, implying each 
$\Gamma_{\epsilon_n, j}$ contains one positive and one negative index. 
On the other hand, 
$\widetilde \mu_0^2<0$ while $\widetilde \mu_j^2 > 0$ for $j>0$. 
Furthermore, we can explicitly compute these indices when $j=1,2,
\dots, n$, in which case $k=1$ and $\lambda_j = n-1$ and 
$$\mu_j ^2 = \frac{n^2+2}{2}, \qquad \widetilde \mu_j^2 = 1.$$
Moreover, the fundamental period of the Jacobi field associated to 
$\widetilde \mu_0$ is 
\begin{equation} \label{cyl_period}
T_{cyl} = T_{\epsilon_n} = \frac{2\pi}{\widetilde \mu_0} , \qquad 
\widetilde \mu_0 = \frac{1}{2} \sqrt{ \sqrt{n^4 - 64n +64}
- n(n-4)+8}.
\end{equation}
To summarize, we have proved the following lemma. 
\begin{lemma} 
We have  
\begin{equation} \label{cyl_indicial_roots1}
\Gamma_{\epsilon_n,0} = \left \{ 0, \pm \frac{1}{2} 
\sqrt{n(n-4)+8 + \sqrt{n^4-64n+64}} \right \} , 
\end{equation} 
and for $j > 0$ 
\begin{eqnarray} \label{cyl_indicial_roots2}
\Gamma_{\epsilon_n, j} & = & \left \{ \pm \frac{1}{2} \sqrt{(n+2(k-1))^2+4 -
\sqrt{n^4+64(k-1)(n+k-1)}}, \right . \\ \nonumber 
&& \left.  \pm \frac{1}{2} \sqrt{ (n+2(k-1))^2+4 
+\sqrt{n^4+64(k-1)(n+k-1)}}\right \}, 
\end{eqnarray} 
where $k$ is the positive integer corresponding to $\lambda_j = 
k(n-2+k)$. In particular, 
\begin{equation} \label{cyl_indicial_roots3} 
\Gamma_{\epsilon_n, 1} = \cdots = \Gamma_{\epsilon_n, n} 
= \left \{ \pm 1, \pm \sqrt{\frac{n^2+2}{2}} \right \}. 
\end{equation} 
\end{lemma}
\begin{rmk} It follows from \eqref{cyl_fourier_exp1} and 
\eqref{cyl_fourier_exp2} that as $k \nearrow \infty$ we have 
$\gamma_{\epsilon_n, j} \simeq \sqrt{2} (k-1) + 2\sqrt{k-1}$ 
and $\widetilde \gamma_{\epsilon_n, j} \simeq \sqrt{2} (k-1) 
- 2\sqrt{k-1}$. However, we will not need this information later. 
\end{rmk} 

For future calculations we will write the 
set of indicial roots as 
$$\Gamma_\epsilon = \{\dots, -\gamma_{\epsilon,2}, -\gamma_{\epsilon,1}=-1, 
0, \gamma_{\epsilon,1}=1, \gamma_{\epsilon_2}, \dots \}$$
where $\gamma_{\epsilon,j} < \gamma_{\epsilon, j+1} \rightarrow \infty$. We 
will justify later the fact that $\Gamma_\epsilon$ has no 
accumulation points. 
 
\subsection{The Fourier-Laplace transform} 

The following transform, defined in \cite{MPU}, plays a key role in our understanding of the 
mapping properties of $L_\epsilon$ and $\{ L_{\epsilon,j} \}$. 

\begin {defn}
Let $\gamma \in \R$ and let 
$w \in W^{k,2}_\gamma((0,\infty) \times \Ss^{n-1})$. Extend $w$ 
to be $0$ in the half-space $\{ t< 0 \}$ and define 
\begin{equation} \label{fourier_laplace_defn}
\mathcal{F}_\epsilon(w)(t,\xi,\theta) = \widehat w(t, \xi, \theta) = \sum_{k=-\infty}^\infty
e^{-i\xi k} w  (t+k T_\epsilon, \theta ).
\end{equation}
Here $\xi \in \{ \eta +i\nu \in \C : \nu <-\gamma T_\epsilon\}$. 
\end{defn} 

\begin{lemma} \label{fourier_laplace_conv}
The sum in \eqref{fourier_laplace_defn} converges uniformly and 
absolutely when $w \in 
W^{k,2}_\gamma((0,\infty) \times \Ss^{n-1})$ and $\nu = \Im(\xi) < 
-\gamma T_\epsilon$. 
Equivalently, 
$$ w \in W^{k,2}_\gamma ((0,\infty) \times \Ss^{n-1}) \Rightarrow 
\mathcal{F}_\epsilon(w) \in \mathcal{C}^\omega (
\{ \Im(\xi) <-\gamma T_\epsilon \}, W^{k,2}
((0,\infty) \times \Ss^{n-1}).$$
\end {lemma}

\begin{proof} 
We have seen that $w \in W^{k,2}_\gamma ((0,\infty) \times \Ss^{n-1})$ 
implies $|w(t,\theta)| = \mathcal{O} (e^{\gamma t})$. 
Writing $\xi = \eta + i\nu$, with $\eta, \nu \in \R$, we have 
\begin{eqnarray*}
| \mathcal{F}_\epsilon (w) (t,\xi,\theta) | & \leq &
\sum_{k=-\infty}^\infty  \left | e^{-i(\eta + i\nu) k} w(t+ kT_\epsilon ,\theta)
\right | = \sum_{k=-\infty}^{\infty} e^{k \nu}|w(t+kT_\epsilon,\theta)| \\
& \leq & C e^{\gamma t} \sum_{k=-\infty}^{\infty} e^{k(\nu +\gamma T_\epsilon)}.
\end {eqnarray*}
First observe that, since we have extended $w$ to be $0$ in the region $\{ t< 0\}$, 
each choice of $\xi = \eta + i \nu$ only gives finitely many nonzero terms 
with $k<0$, and so we only must resolve the convergence when $k \rightarrow 
\infty$. In this case, all exponents are negative precisely when $\nu < 
-\gamma T_\epsilon$. 
\end {proof}

Heuristically, the parameter $\xi$ (more specifically $\nu = \Im(\xi)$) 
allows us to move the weight as a parameter in the  function space to 
one in the operator. 

One can invert this transform, but (as expected) one must choose a 
branch in of the inversion. 
\begin {lemma} Let $w \in W^{k,2}_\gamma((0,\infty) \times \Ss^{n-1})$ and let 
$\nu < -\gamma T_\epsilon$. For each $t$ choose $l \in \mathbf{Z}$ and 
$\tilde t \in [0,T_\epsilon)$ so that $t = \tilde t + l T_\epsilon$. Then 
\begin {equation} \label{fourier_laplace_inv}
w(t,\theta) = \frac{1}{2\pi} \int_{\eta = 0}^{2\pi} 
e^{ilT_\epsilon (\eta + i\nu)}
\widehat w(\tilde t, \eta + i\nu, \theta) d\eta .
\end {equation} 
\end {lemma} 

\begin {proof} Writing $\xi = \eta + i\nu$ we have  
\begin {eqnarray*} 
\frac{1}{2\pi} \int_{\eta = 0}^{2\pi} 
e^{il \xi} \widehat w(\tilde t, \xi, \theta) d\eta 
& = & \frac{1}{2\pi} \int_{\eta = 0}^{2\pi} 
e^{il \xi} \sum_{k=-\infty}^\infty 
e^{-ik \xi} w(\tilde t + k T_\epsilon, \theta) d\eta \\ 
& = & \sum_{k=-\infty}^\infty \frac{1}{2\pi} \int_{\eta = 0}^{2\pi}
e^{i(\eta + i\nu)(l-k)} w(\tilde t + kT_\epsilon, \theta) d\eta \\ 
& = & \frac{w(\tilde t + l T_\epsilon, \theta)} {2\pi} 
\int_{\eta = 0}^{2\pi} d\eta  \\ 
& = & w(t, \theta). \end {eqnarray*}
Here we have used the fact that $\nu < -\gamma T_\epsilon$ to 
allow us to interchange the sum and the integral. 
\end {proof} 
In fact, we can 
treat $\nu$ as a parameter in this inversion, and we see that 
changing $\nu$ alters the weight of the transformed function. 
We make this explicit with a version of the Parseval-Plancherel
identity. 
\begin{lemma} 
For each $\theta \in \Ss^{n-1}$ and $\nu \in \R$ we have 
\begin{equation} \label{parseval1}
\| \widehat w(\cdot, \cdot + i\nu, \theta) \|^2_{L^2([0,T_\epsilon]
\times [0,2\pi])} \simeq 2\pi \| w(\cdot, \theta) \|_{L^2_{\nu/T_\epsilon} 
(\R)}, 
\end{equation}
where $a \simeq b$ means both $a = \mathcal{O} ( b)$ and $b = 
\mathcal{O} (a)$.
\end{lemma} 

\begin{proof} We compute 
\begin{eqnarray*} 
\int_0^{T_\epsilon} \int_0^{2\pi} |\widehat w(t,\eta+i\nu,\theta)|^2 d\eta
dt & = & \int_0^{T_\epsilon} \int_0^{2\pi} \left ( \sum_{k=-\infty}^\infty 
e^{-ik\eta} e^{k\nu} w(t+kT_\epsilon,\theta) \right ) \left ( \sum_{l=-\infty}
^\infty e^{il\eta} e^{l\nu} w(t+lT_\epsilon,\theta) \right ) d\eta dt \\ 
& = & \int_0^{T_\epsilon} \int_0^{2\pi} \sum_{k=-\infty}^\infty 
\sum_{l=-k}^k \binom{k}{l} e^{i(l-k)\eta} e^{(k+l)\nu} w(t+kT_\epsilon, \theta) 
w(t+lT_\epsilon,\theta) d\eta dt \\ 
& = & \int_0^{T_\epsilon} \int_0^{2\pi} \sum_{k=-\infty}^\infty 
e^{2\nu k} (w(t+kT_\epsilon, \theta))^2 d\eta dt \\ 
& \simeq & 2\pi \int_\R (e^{\nu t/T_\epsilon} w(t,\theta))^2 dt. 
\end{eqnarray*} 
Observe that the integrals of all the cross-terms in the sum all vanish 
because $\int_0^{2\pi} e^{ikt} dt = 0$ for each $k \in \mathbf{Z} 
\backslash \{ 0 \}$. 
\end{proof} 
Evaluating the computation above with the choice $\nu=0$ 
we find 
\begin{cor} For each $\theta \in \Ss^{n-1}$ we have 
\begin{equation} \label{parseval2}
\|\widehat w (\cdot, \cdot,\theta) \|^2_{L^2([0,T_\epsilon] \times 
[0,2\pi])} = 2\pi \| w \|^2_{L^2(\R)}. 
\end{equation} 
\end{cor}
\begin{proof} Take $\nu=0$ in \eqref{parseval1}. \end{proof}

Furthermore one can reindex the sum in \eqref{fourier_laplace_defn}
to obtain 
\begin {eqnarray} \label{fourier_laplace_period}
\widehat w(t +T_\epsilon, \xi, \theta) & = & \sum_{k=-\infty}^\infty e^{-ik \xi} 
w(t+T_\epsilon + k T_\epsilon, \theta) \\ \nonumber 
& = & \sum_{k=-\infty}^\infty e^{-ik\xi} w(t+ (k+1)T_\epsilon,\theta) \\ \nonumber 
& = & \sum_{l=-\infty}^\infty e^{-i (l-1)\xi} w(t+lT_\epsilon,\theta) \\ \nonumber 
& = & e^{i\xi} \widehat w(t,\xi, \theta) , \end{eqnarray}
which we can write either as $\widehat w(t,\xi,\theta) = e^{-i\xi}
\widehat w(t+T_\epsilon, \xi,\theta)$ or as $w(t+T_\epsilon,\theta) = 
\mathcal{F}_\epsilon^{-1} (e^{i\xi}\mathcal{F}_\epsilon(w)) (t,\theta)$.
In more geometric/invariant language, this last formula states $\widehat w$ is 
a section of the flat bundle $\Ss^1 \times \Ss^{n-1}$ with holonomy $\xi$ 
around the $\Ss^1$ loop. 

\begin{cor} The Fourier-Laplace transform gives a direct 
integral decomposition 
\begin{equation} \label{integral_decomposition} 
L^2 (\R \times \Ss^{n-1} ) = \int_{\eta \in [0,2\pi]}^\oplus 
L^2_\eta ([0,T_\epsilon] \times \Ss^{n-1}) d\eta, 
\end{equation}
where $L^2_\eta([0,T_\epsilon] \times \Ss^{n-1})$ is the 
$L^2$-completion of 
$$\left \{ w \in \mathcal{C}^0 ([0, T_\epsilon] \times \Ss^{n-1}): 
w(T_\epsilon, \theta) = e^{i T_\epsilon \eta} 
w(0,\theta) \right \}.$$
\end{cor} 
\begin{proof} Combine \eqref{fourier_laplace_period} with 
\eqref{parseval1} and \eqref{fourier_laplace_inv}.
\end{proof}

\subsection{Spectral bands of the Jaocbi operator of the Delaunay 
metrics} 

At this point we use the decomposition \eqref{integral_decomposition} 
to prove a spectral gap result for the Jacobi operator of a 
Delaunay metric. Much of this discussion borrows from \cite{MP_bifur}.

We restrict attention in this section to the space of quasi-periodic 
functions. For each $\eta \in \R$ and nonnegative integer $k$ define 
$W^{k,4}_\eta ([0, T_\epsilon])$ to be the $W^{k,4}$-closure of 
the space of smooth functions on $[0, T_\epsilon]$ subject to the 
boundary conditions 
\begin{equation} \label{qp_bndry_cond}
\frac{d^l w}{dt^l}  (T_\epsilon)  = e^{iT_\epsilon \eta}  
\frac{d^l w}{dt^l} (0), \qquad l=0,1,\dots, k-1 .
\end{equation} 
and denote by $L_{\epsilon, j, \eta}$ the restriction 
$$L_{\epsilon, j, \eta} = L_{\epsilon,j} : W^{4,2}_\eta ([0, T_\epsilon]) 
\rightarrow L^2_\eta ([0,T_\epsilon ]).$$

In order to use the decomposition \eqref{integral_decomposition} we 
define the following twisted operator. 
To begin we define $\widehat L_\epsilon (\xi)$ by 
$\widehat L_\epsilon (\xi) (\widehat v) = \widehat {L_\epsilon (v)}$, 
or $\widehat L_\epsilon (\xi) = \mathcal{F}_\epsilon \circ 
L_\epsilon \circ \mathcal{F}_\epsilon^{-1}$. 
Using \eqref{fourier_laplace_period} we see 
\begin{eqnarray*} 
\widehat L_\epsilon (\xi) (e^{i\xi} \widehat w)
(t,\xi,\theta) & = & \widehat L_\epsilon (\xi) (\widehat w)(t+T_\epsilon, \xi,\theta) 
= \widehat {L_\epsilon (w)} (t+T_\epsilon,\xi,\theta) \\
& = &  e^{i\xi} \widehat {L_\epsilon  v} (t,\xi,\theta)  
= e^{i\xi} \widehat L_\epsilon (\xi) (\widehat v)(t,\xi,\theta) , 
\end{eqnarray*}
which we can rearrange to read 
$$e^{-i \xi} \widehat L_\epsilon (\xi) (e^{i  \xi} \widehat v) 
 = \widehat L_\epsilon (\xi) (\widehat v) .$$ 
This last transformation rule allows us to define the 
twisted operator 
\begin{equation} \label{twisted_op_defn} 
\widetilde L_\epsilon (\xi) (\widehat v) 
= e^{i\xi  t} \mathcal{F}_\epsilon
\circ L_\epsilon \circ \mathcal{F}_\epsilon^{-1} (e^{-i\xi  t} 
\widehat v) ,
\end{equation} 
which is now a well-defined operator 
$$\widetilde L_\epsilon (\xi) : W^{k+4,2} (\Ss^1 \times \Ss^{n-1}) 
\rightarrow W^{k,2} (\Ss^1 \times \Ss^{n-1}),$$
for each value of the paramater $\xi \in \C$. Here we identify $\Ss^1 
= \R / T_\epsilon \mathbf{Z}$. Our key point here is that $\widetilde L_\epsilon 
(\xi)$ act on the {\em same} function space for each value of $\xi$. 

Observe that $\widehat L_\epsilon$ has the same coordinate 
expression as $L_\epsilon$. We can again decompose 
$\widehat L_\epsilon$ and $\widetilde L_\epsilon$ into Fourier 
components, obtaining $\widehat L_{\epsilon,j}$ and 
$\widetilde L_{\epsilon,j}$. In particular, by \eqref{fourier_laplace_period} 
the restriction of 
$\widehat L_{\epsilon,j}(\eta)$ to the interval $[0,T_\epsilon]$ 
is exactly the operator $L_{\epsilon,j,\eta}$ defined above. 

For each $\epsilon$, $j$, and $\eta$ the operator $L_{\epsilon, j, \eta}$ 
is a fourth order ordinary differential operator and we denote its eigenvalues 
by  $\sigma_k(\epsilon,j,\eta)$ for $k=0,1,2,\dots$. Furthermore $L_{\epsilon,j, 0} 
= L_{\epsilon, j, 2\pi}$ 
for each $\epsilon$ and $j$, so we may think of 
$$\sigma_k (\epsilon, j, \cdot) : \Ss^1 \rightarrow \R.$$
We denote the image of this eigenvalue maps by 
\begin{equation} \label{jth_spec_band_defn}
B_k(\epsilon, j) = \{ \sigma \in \R : \sigma = \sigma_k (\epsilon,j,\eta) 
\textrm{ for some }\eta \in [0,2\pi/T_\epsilon] \}
\end{equation} 
as the $j$th {\bf spectral band} of $L_{\epsilon,j}$. 

\begin{lemma} 
Each band $B_k(\epsilon,j)$ is a nondegenerate interval. 
\end{lemma} 

\begin{proof}
Each $L_{\epsilon,j}$ is a fourth order ordinary differential operator, and so 
the ODE $L_{\epsilon,j} v = \sigma v$ has a four-dimensional solution space. If 
the function $\sigma_k(\epsilon,j, \cdot)$ is constant on the interval 
$[0,2\pi]$ then $L_{\epsilon,j} v = \sigma v$ must have an infinite 
dimensional solution 
space for $\sigma \in B_k(\epsilon,j)$. which is impossible. We conclude 
that no band $B_k (\epsilon,j)$ may collapse to a single point. 
\end{proof}

The eigenfunction $w$ corresponding to the eigenvalue $\sigma_k(\epsilon,j,\eta)$
satisfies $w(t+ 2\pi/T_\epsilon) = e^{i\eta} w(t) = e^{(2\pi-\eta) i} w(t)$ and so 
$\bar w(t+2\pi) 
= e^{-i\eta} \bar w(t)$. However, the coefficients of the ordinary differential 
operator $L_{\epsilon,j}$ are real, so 
$$\sigma_k \left (\epsilon, j , 2\pi - \eta \right ) = \sigma_k (\epsilon, 
j, \eta) $$ 
and we may as well restrict $\sigma$ to the half-circle corresponding to 
$0 \leq \eta \leq \pi$. 

It follows from Floquet theory \cite {MW} that the band functions 
$B_{2k}$ are nondecreasing for each $k \in \mathbf{Z}$ while 
$B_{2k+1}$ are all nonincreasing, so that for each $\epsilon$, $j$ and $k$
we have 
\begin{equation} \label{band_ordering1} 
\sigma_0(\epsilon,j,0) \leq \sigma_0(\epsilon,j, \pi) 
\leq \sigma_1 (\epsilon,j,\pi) \leq \sigma_1 (\epsilon,
j,0) \leq \dots.
\end{equation} 
This in turn implies the bands all have the structure 
\begin{eqnarray} \label{band_ordering2} 
B_{2k}(\epsilon,j) (\epsilon,j) & = & [\sigma_{2k} (\epsilon,j,0), \sigma_{2k}
(\epsilon, j, \pi)], \\ \nonumber  
B_{2j+1} (\epsilon,j) & = & [\sigma_{2k+1} (\epsilon,j, \pi), 
\sigma_{2k+1} (\epsilon, j, 0)].
\end{eqnarray} 

We can related to bands $B_k(\epsilon,0)$ to the bands $B_k(\epsilon,j)$ 
using the identity 
\begin{equation} \label{higher_bands1}
L_{\epsilon,j} = L_{\epsilon,0} - 2\lambda_j \frac{d^2}{dt^2} + 
\frac{n(n-4)}{2} \lambda_j + \lambda_j^2. 
\end{equation} 
Let $w$ be an eigenvalue of $L_{\epsilon,j, \eta}$, so that 
\eqref{higher_bands1} implies 
\begin{equation} \label{higher_bands2} 
\sigma_k(\epsilon,j,\eta) w = L_{\epsilon,j} w = L_{\epsilon,0} w 
- 2\lambda_j \ddot w + \frac{n(n-4)}{2} \lambda_j w + \lambda_j^2 w.
\end{equation} 
Writing $w = \sum_{l=0}^\infty \alpha_l w_l$ where 
$L_{\epsilon,0} w_l = \sigma_l(\epsilon,0,\eta) w_l$ we 
rewrite \eqref{higher_bands2} as 
$$
\sum_l \alpha_l \sigma_k(\epsilon,j,\eta) w_l = 
\sum_l \alpha_l\left ( \sigma_l (\epsilon,0,\eta) w_l 
- 2\lambda_j \ddot w_l + \frac{n(n-4)}{2} \lambda_j w_l 
+ \lambda_j^2 w_l\right ),$$
which in turn gives us 
\begin{equation} \label{higher_bands3} 
2\lambda_j \ddot w_l = -\left ( \sigma_k(\epsilon,j,\eta) - 
\sigma_l (\epsilon, 0,\eta) - \frac{n(n-4)}{2} \lambda_j - 
\lambda_j^2 \right )w_l.
\end{equation} 
This last eigenvalue equation admits quasi-periodic solutions 
only if 
\begin{equation} \label{higher_bands4} 
\sigma_k(\epsilon,j,\eta)  > \sigma_l (\epsilon,0,\eta) + \frac{n(n-4)}{2}
\lambda_j + \lambda_j^2.\end{equation}
We have just proved the following lemma. 
\begin{lemma} 
For any positive integer $j$ we have the lower bound
\begin{equation} \label{higher_bands5}
\sigma_k(\epsilon,j,0) > \sigma_0 (\epsilon,0,\eta) + \frac{n(n-4)}{2}
\lambda_j + \lambda_j^2 \geq \sigma_0(\epsilon,0,\eta) + \frac{(n-1)}{2} 
(n^2-2n-2)  . 
\end{equation} 
\end{lemma}

Our main characterization of the spectral bands is the following Proposition. 
\begin{prop} 
For each $\epsilon \in (0, \epsilon_n]$ we have 
\begin{equation} \label{bottom_bands1} 
- \frac{n(n^2-4)}{2} \left ( \frac{1}{T_\epsilon} \int_0
^{T_\epsilon} v_\epsilon^{\frac{2n}{n-4}} dt \right )^{4/n}  \leq 
\sigma_0(\epsilon,0,0) < 0
\end{equation} 
and 
\begin{equation} \label{bottom_bands2} 
\textrm{either }\sigma_1(\epsilon,0,0) = 0 \textrm{ or }
\sigma_2(\epsilon,0,0) = 0.
\end{equation} 
\end{prop} 

\begin{proof} Observe that $\dot v_\epsilon$ is a periodic 
solution of the ODE $L_{\epsilon,0} (\dot v_\epsilon) = 0$, so 
it must be an eigenfunction with associated eigenvalue $0$, subject 
to periodic boundary conditions, {\it i.e.} $\eta = 0$. This 
eigenfunction has precisely two modal domains within the interval $[0,
T_\epsilon]$, so it must correspond either to $\sigma_1(\epsilon,0)$ 
or to $\sigma_2(\epsilon, 0)$. We don't have enough information at this 
point to distinguish these two cases. 

The function $v_\epsilon$ is also $T_\epsilon$-periodic, and so 
is an appropriate test function for $\sigma_0(\epsilon,0,0)$. We have 
\begin{eqnarray*} 
L_{\epsilon,0} (v_\epsilon) & = & \ddddot{v_\epsilon} - \left ( \frac{n(n-4)+8}{2}
\right ) \ddot v_\epsilon + \frac{n^2(n-4)^2}{16} v_\epsilon - 
\frac{n(n+4)(n^2-4)}{16} v^{\frac{n+4}{n-4}} \\ 
& = & \ddddot v_\epsilon + \left ( \frac{n(n-4)+8}{2} \right ) \ddot v_\epsilon 
+ \frac{n^2(n-4)^2}{16} v_\epsilon - \frac{n(n-4)(n^2-4)}{16} 
v_\epsilon^{\frac{n+4}{n-4}} - \frac{n(n^2-4)}{2} v_\epsilon^{\frac{n+4}{n-4}} \\ 
& = & -\frac{n(n^2-4)}{2} v_\epsilon^{\frac{n+4}{n-4}} < 0 ,
\end{eqnarray*} 
and so 
$$ \sigma_0(\epsilon, 0,0) \leq \frac{\int_0^{T_\epsilon} v_\epsilon L_{\epsilon,0} 
(v_\epsilon) dt}{\int_0^{T_\epsilon} v_\epsilon^2 dt} = -\frac{n(n^2-4)}{2} 
\frac{\int_0^{T_\epsilon} v_\epsilon^{\frac{2n}{n-4}} dt}{\int_0^{T_\epsilon} 
v_\epsilon^2 dt} < 0,$$
which gives the upper bound in \eqref{bottom_bands1}. 

On the other hand, combining the uniqueness theorem of \cite{FK} 
and the variational characterization of the Deleaunay solution in Section 5 
of \cite{JX} 
we we that (up to translations) $v_\epsilon$ is the unique minimizer of the 
functional 
$$W^{4,2}_0([0,T_\epsilon]) \ni v \mapsto \frac{\int_0^{T_\epsilon} 
\ddot v^2 + \left ( \frac{n(n-4)+8}{2} \right ) \dot v^2 + \frac{n^2(n-4)^2}{16} 
v^2 dt }{\left ( \int_0^{T_\epsilon} v^{\frac{2n}{n-4}} dt \right )^{\frac{n-4}{n}}} . 
$$ 
By \eqref{paneitz_ode1} we then have 
\begin{equation} \label{bottom_bands3} 
\frac{\int_0^{T_\epsilon} \ddot v^2 + \left (\frac{n(n-4)+8}{2} 
\right )\dot v^2 + \frac{n^2(n-4)^2}{16} v^2 dt} {\left ( 
\int_0^{T_\epsilon} v^{\frac{2n}{n-4}} dt \right )^{\frac{n-4}{n}}} 
\geq \frac{n(n-4)(n^2-4)}{16} \left ( \int_0^{T_\epsilon} 
v_\epsilon^{\frac{2n}{n-4}} dt \right )^{4/n} 
\end{equation} 
for each $v \in W^{4,2}_0([0,T_\epsilon])$, which then gives 
\begin{eqnarray} \label{bottom_bands4} 
\int_0^{T_\epsilon} v L_{\epsilon,0} (v) dt & = & \int_0^{T_\epsilon} 
\ddot v^2 + \left ( \frac{n(n-4)+8}{2} \right )\dot v^2 + \frac{n^2(n-4)^2}
{16} v^2 - \frac{n(n+4)(n^2-4)}{16} v^{\frac{2n}{n-4}} dt \\ \nonumber 
& \geq & \frac{n(n-4)(n^2-4)}{16} \left ( \int_0^{T_\epsilon} v_\epsilon^{\frac{2n}{n-4}}
dt \right )^{4/n} \left ( \int_0^{T_\epsilon} v^{\frac{2n}{n-4}} dt 
\right )^{\frac{n-4}{n}} - \frac{n(n+4)(n^2-4)}{16} \int_0^{T_\epsilon} 
v^{\frac{2n}{n-4}} dt \\ \nonumber 
& = & \left ( \int_0^{T_\epsilon} v^{\frac{2n}{n-4}} dt \right )^{\frac{n-4}{n}}
\left ( \frac{n(n-4)(n^2-4)}{16}  \left ( \int_0^{T_\epsilon} v_\epsilon
^{\frac{2n}{n-4}} dt \right )^{4/n} - \frac{n(n+4)(n^2-4)}{16} \left ( 
\int_0^{T_\epsilon} v^{\frac{2n}{n-4}} dt \right )^{4/n} \right ) .
\end{eqnarray}
H\"older's inequality with exponents $\frac{n}{n-4}$ and $n/4$ implies 
$$\int_0^{T_\epsilon} v^2 dt \leq T_{\epsilon}^{4/n} \left ( \int_0^{T_\epsilon} 
v^{\frac{2n}{n-4}} dt \right )^{\frac{n-4}{n}},$$
which we combine with \eqref{bottom_bands4} to see 
\begin{eqnarray} \label{bottom_bands5} 
\int_0^{T_\epsilon} v L_{\epsilon,0} (v) dt & \geq & 
\frac{n(n-4)(n^2-4)}{16}  T_\epsilon^{-4/n} \int_0^{T_\epsilon} 
v^2 dt 
\left ( \int_0^{T_\epsilon} v_\epsilon
^{\frac{2n}{n-4}} dt \right )^{4/n}  \\ \nonumber 
&& - \frac{n(n+4)(n^2-4)}{16} T_\epsilon^{-4/n} \int_0^{T_\epsilon}
v^2 dt  \left ( 
\int_0^{T_\epsilon} v^{\frac{2n}{n-4}} dt \right )^{4/n}  
\end{eqnarray}
for each $v \in W^{4,2}_0([0,T_\epsilon])$. Finally, we are free 
to choose a scale for our test function $v$, and we normalize so 
that $\int_0^{T_\epsilon} v^{\frac{2n}{n-4}} dt = \int_0^{T_\epsilon} 
v_\epsilon^{\frac{2n}{n-4}} dt$. Using this normalization 
\eqref{bottom_bands5} becomes 
\begin{eqnarray*}
 \int_0^{T_\epsilon} v L_{\epsilon,0}(v) dt & \geq & \left ( 
 \frac{n(n-4)(n^2-4)}{16} - \frac{n(n+4)(n^2-4)}{16} \right ) 
 T_\epsilon^{-4/n} \left ( \int_0^{T_\epsilon} v^2 dt \right ) 
 \left ( \int_0^{T_\epsilon} 
 v_\epsilon^{\frac{2n}{n-4}} dt \right )^{4/n} \\ 
 & = & -\frac{n(n^2-4)}{2} T_\epsilon^{-4/n} \left ( \int_0^{T_\epsilon} 
 v^2 dt \right ) \left ( \int_0^{T_\epsilon} 
 v_\epsilon^{\frac{2n}{n-4}} dt \right )^{4/n} ,
 \end{eqnarray*} 
 which gives the lower bound in \eqref{bottom_bands1}.
\end{proof}

\begin{cor} 
We have $B_k(\epsilon,0) \subset (0,\infty)$ for each $k \geq 3$ and 
$B_k(\epsilon,0) \subset [0,\infty)$ for each $k \geq 2$. 
\end{cor}
\begin{proof} This follows from the previous proposition 
and \eqref{band_ordering2}. 
\end{proof}

\begin{cor} For each positive integer $j$ we have $B_k(\epsilon,j) \subset 
(0,\infty)$. 
\end{cor} 

\begin{proof} 
When $j >n$ we have $\lambda_j \geq 2n$, and in this 
case \eqref{higher_bands5} gives 
$$\sigma_k (\epsilon,j,0) > \sigma_0(\epsilon,0,0) + n^3$$
for all $k$. 
On the other hand, $0<v_\epsilon <1$, so the lower bound 
in \eqref{bottom_bands1} gives us 
$$\sigma_0 (\epsilon,0,0) \geq -\frac{n(n^2-4)}{2} \left ( \frac{1}{T_\epsilon}
\int_0^{T_\epsilon} v_\epsilon^{\frac{2n}{n-4}} dt \right )^{4/n} \geq 
-\frac{n(n^2-4)}{2} ,$$
which implies 
$$\sigma_k(\epsilon,j,0) > \sigma_0(\epsilon,0,0) + n^3 \geq n^3-\frac{n(n^2-4)}{2} 
> 0.$$
and the corollary follows. 

The case of $1 \leq j \leq n$ requires some more attention. Luckily, 
we know from \eqref{low_paneitz_mode2} and \eqref{low_paneitz_mode3}
that 
$$ L_{\epsilon,j} (w_j^\pm) = 0, \qquad 
w_j^\pm = e^{\pm t} \left ( \pm \dot v_\epsilon + \frac{n-4}{2}
v_\epsilon \right ) + \mathcal{R}_\pm $$
where $\mathcal{R}_+ = \mathcal{O}(1)$  
and $\mathcal{R}_- = \mathcal{O}(e^{-2t})$. By 
\eqref{trans_nodal_domains2} these
are positive, periodic solutions, and so must 
correspond to the bottom of the spectrum of $L_{\epsilon, j, \eta}$. 
\end{proof} 

The following lemma relates the spectral bands $B_k(\epsilon, j)$ and 
the set of indicial roots $\Gamma_{\epsilon,j}$ of the 
operator $L_{\epsilon,j}$. 
\begin{lemma} 
We have $0 \in B_k(\epsilon,j)$ for some $k$ if and only if 
the ODE $L_{\epsilon,j} w = 0$ admits a quasi-periodic solution.  
\end{lemma} 

\begin{proof} 
Let $\widehat v$ satisfy 
$$0 = L_{\epsilon,j,\eta} \widehat v = e^{i\eta t}\mathcal{F}_\epsilon( L_{\epsilon,j} 
(\mathcal{F}_\epsilon^{-1} (e^{-i\eta t} \widehat v)))\Rightarrow 
L_{\epsilon,j} (\mathcal{F}_\epsilon^{-1} (e^{-i\eta t} \widehat v)) = 0. $$
In particular, $v = \mathcal{F}_\epsilon^{-1} (e^{-i\eta t} \widehat v)$ 
satisfies $L_{\epsilon,j} (v) = 0$. In addition, \eqref{fourier_laplace_period}
implies $v(t+T_\epsilon,\theta) = e^{i\eta} v(t,\theta)$, and so 
$v$ must be quasi-periodic. 
\end{proof} 

We summarize the important conclusions of this section with 
the following Corollary. 
\begin{cor} \label{low_indicial_roots_cor}
The indicial root $0$ is isolated, {\it i.e.} there exists $\delta>0$ 
such that no other indicial roots lie in the interval $(-\delta, \delta)$. 
Moreover, any Jacobi field with sub-exponential growth ({\it i.e.} 
tempered) must be a linear combination of $w_0^+$ and $w_0^-$, 
the Jacobi fields generated by translations along the axis and 
changes of the Delaunay parameter. 
\end{cor}

\subsection{Mapping properties of the linearized operator}

We have already seen that for each $\delta \in \R$ and $k \in \mathbf{N}$ 
the mapping 
$$L_\epsilon: W^{k+4,2}_\delta ((0, \infty) \times \Ss^{n-1}) \rightarrow 
W^{k,2}_\delta ((0, \infty) \times \Ss^{n-1}) $$ 
is a linear, elliptic operator and it has bounded coefficients. 

\begin{prop} \label{nondegeneracy} 
The operator 
$$L_\epsilon : W^{k+4,2}_\delta ((0,\infty) \times \Ss^{n-1}) \rightarrow 
W^{k,2}_\delta ((0,\infty) \times \Ss^{n-1})$$ 
is Fredholm provided $\delta \not \in \Gamma_\epsilon$, where 
$\Gamma_\epsilon$ is given in \eqref{indicial_defn}. 
\end{prop} 

Our proof follows that of Proposition 4.8 of \cite{MPU}. 
\begin{proof} 
We recall the twisted operator define in \eqref{twisted_op_defn}  
depending on a parameter $\xi \in \C$, namely 
$$\widetilde {L}_\epsilon (\xi) : W^{k+4,2} (\Ss^1 \times \Ss^{n-1}) 
\rightarrow W^{k,2} (\Ss^1 \times \Ss^{n-1}), \qquad 
\widetilde{L}_\epsilon (\xi) (\widehat v) = e^{i\xi t} \mathcal{F}_\epsilon 
\circ L_\epsilon \circ \mathcal{F}_\epsilon^{-1} (e^{-i\xi t} \widehat v).$$
We then use the analytic Fredholm theorem to prove the twisted 
operator is Fredholm away from a discrete set of poles in 
$\C$. Unwinding definition, we then translate the Fredholm 
property of the twisted operator to corresponding properties of 
$L_\epsilon$ on a weighted function space. 

The operator $\widetilde L_{\epsilon}(\xi)$ is linear, bounded, 
and elliptic for each choice of $\xi$, and depends on $\xi$ holomorphically. 
Thus, by the analytic 
Fredholm theorem (see Section 5.3 of \cite{Mel}) $\widetilde L_{\epsilon}(\xi)$ is either 
never Fredholm for any value $\xi$ or it is Fredholm for $\xi$ not in a certain 
discrete set. We choose $\xi = \eta \in (0,2\pi)$ and suppose there exists 
$\widehat v$ such that $\widetilde L_\epsilon (\eta) (\widehat v) =0$, then 
$L_\epsilon (v) = 0$ where $v = \mathcal{F}_\epsilon ^{-1}
(e^{-i\eta t} \widehat v)$. Then $v$ is quasi-periodic, in that 
$v(t+T_\epsilon, \theta) = e^{i\eta } v(t,\theta)$, and 
so in particular $v$ is bounded. However, by 
Corollary \ref{low_indicial_roots_cor} any bounded Jacobi field must be 
a multiple of $w_0^+$, which is not quasi-periodic. 
Thus $\widetilde L_\epsilon (\eta)$ 
is injective. However, this operator is formally self-adjoint, so it is also 
surjective, and in particular $\widetilde L_\epsilon (\eta)$ is Fredholm. 

Therefore we may safely apply the analytic Fredholm theorem to conclude 
there exists a discrete set $\mathcal{P} \subset \C$ and a meromorphic operator 
$$\widetilde G_\epsilon (\xi) : W^{k,2} (\Ss^1 \times \Ss^{n-1}) 
\rightarrow W^{k+4,2} (\Ss^1 \times \Ss^{n-1}) $$
such that $\widehat v = \widetilde G_\epsilon (\xi) \circ \widetilde 
L_\epsilon (\xi)(\widehat v)$ so long as $\xi \not \in \mathcal{P}$. 
Moreover, by our construction 
$$\Gamma_\epsilon = \{ \nu \in \R : \nu = \Im (\xi) \textrm{ for 
some } \xi \in \mathcal{P} \}.$$

We unravel this relation to find the Greens operator $G_\epsilon$. Again 
we let $v = \mathcal{F}_\epsilon^{-1}(e^{-i\xi t} \widehat v)$, so 
that 
\begin{eqnarray*}
e^{i\xi t} \mathcal{F}_\epsilon (v) & = & \widehat v   
= \widetilde G_\epsilon (\widetilde L_\epsilon (\widehat v)) \\ 
& = &  \widetilde G_\epsilon (e^{i\xi t} \mathcal{F}_\epsilon 
(L_\epsilon (\mathcal{F}_\epsilon^{-1} ( e^{-i\xi t} 
e^{i \xi t} \mathcal{F}_\epsilon (v))))) \\ 
& = & \widetilde G_\epsilon ( e^{i \xi t} \mathcal{F}_\epsilon 
(L_\epsilon (v))). 
\end{eqnarray*} 
We thus conclude 
\begin{equation} \label{greens_op}
G_\epsilon (\phi) = \mathcal{F}_\epsilon^{-1} (e^{-i\xi T_\epsilon t}
(\widetilde G_\epsilon (e^{i\xi T_\epsilon t} (\mathcal{F}_\epsilon (\phi)))).
\end{equation} 
By our construction we have 
$v = G_\epsilon (\phi) \in W^{k+4,2}_{-\Im(\xi)}
((0,\infty) \times \Ss^{n-1})$, so that 
$$\delta = \Im(\xi) \not \in \Gamma_\epsilon \Rightarrow 
\textrm{ there exist a Greens operator } G_\epsilon : 
W^{k,2}_\delta ((0,\infty) \times \Ss^{n-1}) \rightarrow W^{k+4,2}_\delta 
((0,\infty) \times \Ss^{n-1}),$$
and by the Fredholm alternative we have completed our proof.
\end{proof}

\begin{lemma} The set of indicial roots $\Gamma_\epsilon$ does 
not have any accumulation points. \end{lemma} 

\begin{proof} Recall that each $\gamma \in \Gamma_\epsilon$ is 
the imaginary part of a pole of $\widetilde G_\epsilon$, which 
forms a discrete set 
in $\C$. Furthermore, the operator $\widetilde L_\epsilon(\xi)$ is
unitarily equivalent to $\widetilde L_\epsilon (\xi + 2\pi k)$ for 
each $k \in \mathbf{Z}$, and 
so $\xi$ is a pole of $\widetilde G_\epsilon$ if and only if 
$\xi + 2\pi k$ is as well for each $k \in \mathbf{Z}$.  Thus $\widetilde 
G_\epsilon$ can only have finitely many poles in each horizontal strip. 
\end{proof}

\begin{cor} 
For each $\delta \in (0,1)$ the operator 
$$L_\epsilon : W^{k+4,2}_{-\delta} ((0,\infty) \times \Ss^{n-1}) \rightarrow 
W^{k,2}_{-\delta} ((0,\infty) \times \Ss^{n-1}) $$
is injective and 
$$L_\epsilon : W^{k+4,2}_\delta ((0,\infty) \times \Ss^{n-1}) \rightarrow 
W^{k,2}_\delta ((0,\infty) \times \Ss^{n-1}) $$ 
is surjective. \end{cor} 

\begin{proof} 
In the proof of Proposition \ref{nondegeneracy} we showed that 
$\widetilde L_\epsilon (\xi)$ is injective for each $\xi$ with 
$-1 < \Im(\xi) < 0$,which in turn implies 
$L_\epsilon : W^{k+4,2}_{-\delta} ((0,\infty) \times \Ss^{n-1}) 
\rightarrow W^{k,2}_{-\delta} ((0,\infty) \times \Ss^{n-1})$ 
is injective. The surjectivity statement now follows from 
the injectivity statement and duality because $L_\epsilon: 
W^{k+4,2}_0 ((0,\infty) \times \Ss^{n-1}) \rightarrow W^{k,2}_0
((0,\infty) \times \Ss^{n-1})$ is formally self-adjoint. 
\end{proof} 

We can extract more information from the Fourier-Laplace transform by examining 
how it behaves when the contour we are using to define its inverse crosses a 
pole of $\widetilde G_\epsilon$. We use Proposition 4.14 of \cite{MPU} as a 
model for the following. 

\begin {prop} \label{lin_soln_asymp_expansion}
Let $\phi \in \mathcal{C}^\infty_0 ((0,\infty) \times \Ss^{n-1})$ and let 
$\delta \in (0,1)$ and 
$v \in W^{4,2}_{-\delta} ((0,\infty) \times \Ss^{n-1})$ 
satisfy $L_\epsilon (v) = \phi$. Then $v$ has an asymptotic 
expansion $v = \sum_{k=1}^\infty v_j$ with $L_\epsilon (v_j) = 0$ 
and $v_j \in W^{4,2}_{-\nu}$ for any $\nu < \gamma_{\epsilon,j}$.
\end {prop}

\begin {proof} We begin by transforming the equation. Choose 
$\xi \in \C$ with $\Im (\xi) > \delta>0$ and let 
$$\widetilde v = e^{i\xi  t} \widehat v, \qquad 
\widetilde \phi = e^{i \xi t} \widehat \phi,$$
so that $\widetilde L_\epsilon (\xi) (\widetilde v) = 
\widetilde \phi$. Applying the Greens operator 
$\widetilde G_\epsilon (\xi)$ we have $\widetilde v = 
\widetilde G_\epsilon (\xi)(\widetilde \phi)$. That 
$\phi \in \mathcal{C}_0^\infty ((0,\infty) \times \Ss^{n-1})$ 
implies $\widetilde \phi$ is entire in the $\xi$ variable 
and smooth in the $(t,\theta)$ variables. On the other 
hand, in the half-plane $\Im (\xi) > \delta$ the poles of 
$\widetilde G_\epsilon$ occur at the same points as the 
zeroes of $\widetilde \phi$ (with the same degrees), so 
$\widetilde v$ is analytic in this half-plane. In fact, because 
$\gamma_{\epsilon,1} = 1$, for any 
$\delta' \in (\delta, 1)$ the operator
$\widetilde G_\epsilon$ has no poles in the strip 
$\delta' < \Im (\xi) < \delta$, we can shift the contour integral 
defining $\mathcal{F}_\epsilon^{-1}$ up to $\Im(\xi) = \delta'$ 
for any $\delta' < 1$, and so we may take 
$v \in W^{k+4,2}_{-\nu}((0,\infty) \times \Ss^{n-1})$ for any $\nu < 
1$.

\begin{figure} [h] 
\centering
\begin{tikzpicture}
\coordinate (0) at (2,2); 
\draw[->] (0,3) -- (8,3) coordinate[label = {below:$x$}] (xmax);
\draw[->] (4,-2) -- (4,4) coordinate[label = {right:$y$}] (ymax);
\draw[thick, blue,->>] (2,0.5) -- (4,0.5); \draw[thick,blue] (4, 0.5) -- (6, 0.5); 
\draw[thick, red](2,1.5) -- (4,1.5); \draw[thick, red,<<-] (4,1.5) -- (6,1.5); 
\draw[thick, green,->>] (6,0.5) -- (6, 1.0); \draw[thick,green] (6,1.0) -- (6,1.5);  
\draw[thick,green] (2,0.5) -- (2,1.0); \draw[thick, green,<<-] ( 2,1.0) -- (2,1.5);   
\node at (4,1) {$*$};  
\node at (5.6,-0.8) {pole at $-i$}; 
\draw [->] (5.5,-0.5) -- (4.05,1); 
\end{tikzpicture} 
\caption{This figure the contour we integrate along, surrounding 
the pole corresponding to the first indicial root $1$. 
The horizonal segments have height $-\delta''$ and $-\delta$, and extend 
horizontally from $-\pi$ to $\pi$. The contour integrals along the two 
vertical segments cancel out by the periodicity 
in \eqref{fourier_laplace_inv}.} \label{contour_fig}
\end{figure}
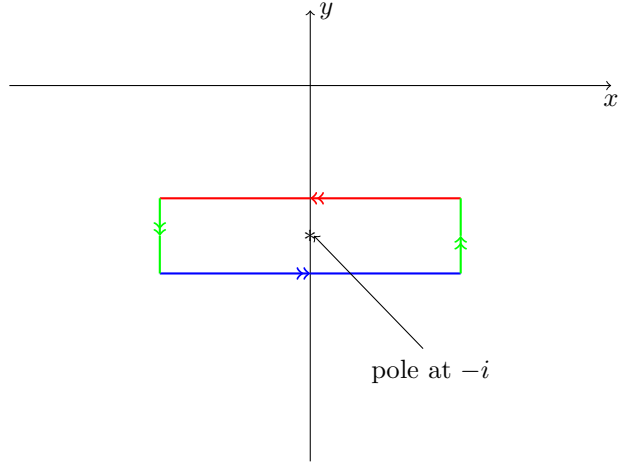

We complete the proof by shifting the contour integral across 
a pole of $\widetilde G_\epsilon$. We sketch this 
contour in Figure \ref{contour_fig}. Choose $\delta'' \in 
(1, \gamma_{\epsilon,2})$ and $\xi''$ such that 
$\Im(\xi'') = \delta''$, and let $\widetilde v'' = 
\widetilde G_\epsilon (\xi'')$. Applying $\mathcal{F}_\epsilon^{-1}$ 
along the two contours $\Im(\xi) = \delta$ and $\Im(\xi) = \delta''$
we thus have $v = v' + v''$, where $v'' \in W^{k+4,2}_{-\delta''} 
((0,\infty) \times \Ss^{n-1})$ and $L_\epsilon (v') = 0$ and 
$v' \in W^{k+4,2}_{-\nu}((0,\infty) \times \Ss^{n-1})$ for 
each $\nu \in (0,1)$. In fact the difference 
$\widetilde v - \widetilde v'' = (\widetilde G_\epsilon (\xi) - 
\widetilde G_\epsilon (\xi'')) (\widetilde \phi)$ is the 
residue of a meromorphic function with a pole at height 
$\Im(\xi) = -1$ about a rectangular contour 
of width $2\pi$ and height $\delta'' - \delta$. We do not see 
the result of the contour integral along the vertical sides of 
the rectangle because all our transformed functions are 
periodic in the real direction with period $2\pi$.  
\end {proof} 

The following is a special case of 
Proposition \ref{lin_soln_asymp_expansion}, corresponding 
to moving the contour integral across only the first pole. 
\begin{cor} \label{lin_decomp_lemma}
Let $\delta \in (0,1)$, let $\phi \in \mathcal{C}^\infty ((0, \infty) 
\times \Ss^{n-1}) \cap L^2_{-\delta} ((0,\infty) \times \Ss^{n-1})$ and 
let $v \in W^{4,2}_\delta ((0,\infty) \times \Ss^{n-1})$ such that 
$L_\epsilon (v) = \phi$. Then there exist $w \in W^{4,2}_{-\delta} 
((0,\infty) \times \Ss^{n-1})$ and $ \psi \in \operatorname{Span} 
\{ w_0^+, w_0^-\}$ such that 
$v = w+\psi$. 
\end{cor}

\section{Simple convergence to a radial solution} 
\label{asymp_sec}

In this section we prove an asymptotic estimate. 
\begin {thm} Let $v \in \mathcal{C}^\infty ((0,\infty) \times 
\Ss^{n-1})$ be a positive solution of \eqref{paneitz_pde2} 
which also satisfies \eqref{pos_lap_cyl_coords}. Then 
either $\limsup_{t \rightarrow \infty} v(t,\theta) = 0$ 
or there exists a Delaunay 
 parameter $\epsilon \in (0,v_{cyl}]$ and translation parameter 
 $T \in [0,T_\epsilon )$ and positive 
 constants $C$ and $\alpha$ such that 
\begin {equation} \label{asymp_symm1}
|v(t,\theta) - v_\epsilon (t+T)| \leq C e^{-\alpha t} .
\end {equation} 
\end {thm} 

\begin{rmk} Strictly speaking, the proof below is not entirely 
necessary, as one should be able to show the estimate 
\eqref{asymp_symm1} is equivalent to the simple 
asymptotics Jin and Xiong derive in \cite{JX}. However, 
we believe the proof is different and interesting 
enough to include. 
\end{rmk} 

\begin {proof} The transformed function $v$ 
satisfies \eqref{paneitz_pde2}.
We let $\tau_k \rightarrow \infty$ and define 
$$ v_k: (-\tau_k, \infty) \times \Ss^{n-1} \rightarrow \infty, \qquad 
v_k(t,\theta) = v(t+ \tau_k, \theta).$$
By the bounds \eqref{apriori_bounds2} 
there are constants $0<c_1< c_2$ (which depend on the solutions $v$)
such that $c_1 < v_k(t,\theta) < c_2$ for all $k$. Moreover, elliptic estimates 
imply the sequence $\{ |\nabla v_k |\}$ is also uniformly bounded, so 
in fact the family $\{ v_k\}$ is also equicontinuous. Thus a subsequence 
converges uniformly on compact sets to a solution $\bar v = 
\lim_{l \rightarrow \infty} v_{k_l} \rightarrow \bar v$. Since 
$\tau_k \rightarrow \infty$ we now have a global solution 
$\bar v : \R \times \Ss^{n-1} \rightarrow (0,\infty)$ of the PDE 
\eqref{paneitz_pde2}. By the classification theorem in \cite{FK} we 
must have $\bar v(t,\theta) = v_\epsilon (t+T)$ for some $\epsilon
\in (0,\epsilon_n]$ and $t \in [0,T_\epsilon)$. 

It remains to show that $\epsilon$ and $T$ do not depend on the 
choice of sequence $\tau_k \rightarrow \infty$ or the choice 
of subsequence $v_{k_l}$. 

By \eqref{poho1} for each $k$ we have 
\begin {eqnarray*} 
n \omega_n \mathcal{H}_{\textrm{rad}}(v_\epsilon)  
& = & \lim_{k \rightarrow \infty} \mathcal{H}_{rad} (v_k) \\ 
& = & \lim_{k \rightarrow \infty} \int_{\{ 1 \} \times \Ss^{n-1}} - \frac{\partial v_k}{\partial t} 
\frac{\partial^3 v_k}{\partial t^3}+ \frac{1}{2} \left ( \frac{\partial^2 v_k}
{\partial t^2} \right )^2 + \left ( \frac{n(n-4) + 8}{4} \right ) \left ( \frac 
{\partial v_k}{\partial t} \right )^2 \\ 
&& \quad + \frac{n^2(n-4)^2}{32} v_k^2 + \frac{(n-4)^2(n^2-4)}{32} 
v_k^{\frac{2n}{n-4}} - \frac{1}{2} (\Delta_\theta v_k)^2- \frac{n(n-4)}{4}
\left | \nabla_\theta \frac{\partial v_k}{\partial t} \right |^2 d\sigma \\ 
& = &  \lim_{k \rightarrow \infty} \int_{\{ \tau_k \} \times \Ss^{n-1}} - \frac{\partial v}{\partial t} 
\frac{\partial^3 v}{\partial t^3}+ \frac{1}{2} \left ( \frac{\partial^2 v}
{\partial t^2} \right )^2 + \left ( \frac{n(n-4) + 8}{4} \right ) \left ( \frac 
{\partial v}{\partial t} \right )^2 \\ 
&& \quad + \frac{n^2(n-4)^2}{32} v^2 + \frac{(n-4)^2(n^2-4)}{32} 
v^{\frac{2n}{n-4}} - \frac{1}{2} (\Delta_\theta v)^2- \frac{n(n-4)}{4}
\left | \nabla_\theta \frac{\partial v}{\partial t} \right |^2 d\sigma \\ 
& = &  \int_{\{ 1 \} \times \Ss^{n-1}} - \frac{\partial v}{\partial t} 
\frac{\partial^3 v}{\partial t^3}+ \frac{1}{2} \left ( \frac{\partial^2 v}
{\partial t^2} \right )^2 + \left ( \frac{n(n-4) + 8}{4} \right ) \left ( \frac 
{\partial v}{\partial t} \right )^2 \\ 
&& \quad + \frac{n^2(n-4)^2}{32} v^2 + \frac{(n-4)^2(n^2-4)}{32} 
v^{\frac{2n}{n-4}} - \frac{1}{2} (\Delta_\theta v)^2- \frac{n(n-4)}{4}
\left | \nabla_\theta \frac{\partial v}{\partial t} \right |^2 d\sigma \\ 
& = & \mathcal{H}_{rad} (v), 
\end {eqnarray*}
so $\epsilon$ does not depend on any choices. 

We complete the proof with the help of the following lemmas. 

\begin{lemma} \label{helper_lemma1}
Let $\theta \in \Ss^{n-1}$ and let $v: (0,\infty) 
\times \Ss^{n-1} \rightarrow (0,\infty)$ solve \eqref{paneitz_pde2}. Then 
\begin{equation} \label{ang_decay1}
\lim_{t\rightarrow\infty} \frac{\partial v}{\partial \theta} (t,\cdot) = 0 
\end{equation} 
uniformly. 
\end{lemma} 

\begin{proof} If this were not the case then there exist $(\tau_j, 
\theta_j)\in (0,\infty) \times \Ss^{n-1}$ and $C>0$ such that 
$\tau_j \rightarrow \infty$ and 
$$\left | \frac{\partial v}{\partial \theta} (\tau_j, \theta_j) \right | 
\geq C > 0.$$
Now translate by $-\tau_j$ and rotate by $-\theta_j$ to obtain a 
new sequence of solutions $v_j$, for which 
$$\left | \frac{\partial v_j}{\partial \theta} (0,\theta_0) \right | \geq C
 > 0.$$
However, by the reasoning above the sequence $v_j$ 
must converge to a Delaunay solution, which does not depend on 
$\theta$ at all, contradicting the lower bound displayed above. 
\end{proof} 

\begin{lemma} \label{helper_lemma2}
Let $\theta \in \Ss^{n-1}$ and let $v$ solve 
\eqref{paneitz_pde2}. The Jacobi field $\phi =  \partial_\theta v$. 
decays exponentially in $t$. 
\end{lemma} 

\begin{proof}
Let $\tau_j \rightarrow \infty$ and let 
$$v_j (t,\theta) = v(t+\tau_j, \theta), \quad A_j = \sup \left \{ \left | 
\frac{\partial v_j}{\partial \theta} \right | : t>0 \right \}.$$
Further suppose there exist $(s_j, \theta_j)$ such that $|\partial_\theta
v_j (s_j,\theta_j)| = |\partial_\theta v(s_j + \tau_j,\theta_j)| = A_j$.  If 
the  sequence $\{ s_j\}$ is unbounded then we can translate further 
to obtain $\phi_j (t,\theta) = \partial_\theta v_j(t+s_j, \theta) = 
\partial_\theta v(t+s_j+\tau_j,\theta)$. Letting $j \rightarrow \infty$ 
we obtain $\phi: \R \times \Ss^{n-1} \rightarrow \R$, which solves 
the linearized equation \eqref{linearization2} about $v_\epsilon (\cdot + T)$. 
By our construction $\phi$ is not identically zero and bounded. However, 
by our construction $\phi \not \in \operatorname{Span}\{ w_0^+ = 
\dot v_\epsilon, w_0^- = \frac{d}{d\epsilon} v_\epsilon \}$, so by 
Corollary \ref{low_indicial_roots_cor} 
$\phi$ must grow exponentially either as $t \rightarrow \infty$ or 
$t \rightarrow -\infty$. We have already shown $\phi(t,\theta)$ must decay 
at some rate as $t\rightarrow \infty$, so there must exist $c>0$ such that  
$|\phi(t,\theta)| \leq c e^{- t}$. 
\end{proof} 

A consequence of Lemma \ref{helper_lemma2} is that 
any choice of $s_j>0$ such that 
$$v_j(s_j, \theta_j) = A_j = \sup \left \{ \left | \frac{\partial v_j}
{\partial \theta} \right | : t>0\right \}$$ 
remains bounded. In particular, there exists a positive integer 
$N$ such that $N T_\epsilon > \sup_{j\in \mathbf{N}} s_j$, where $T_\epsilon$ 
is the period of $v_\epsilon$. We now define the intervals 
$I_N = [0,NT_\epsilon]$ and $J_N = [NT_\epsilon, 2NT_\epsilon]$, 
and observe that the limit Jacobi field $\phi$ we have just constructed 
is bounded for $t>0$ and attains its supremum in $I_N$.

\begin{lemma} \label{helper_lemma3} 
Let $\phi$ be the Jacobi field constructed in Lemma \ref{helper_lemma2}. 
Then there exists $c$ independent of all choices such that $|\phi(t,\theta)| 
\leq c e^{- t}$. 
\end{lemma} 

\begin{proof} 
Expand the Jacobi field $\phi$ we have just constructed 
in Fourier modes. Let $\{ w_j^{\pm}, \widetilde{w}_j^{\pm} \}$ span the
solutions space of \eqref{paneitz_ode3} and let 
$$E' = \operatorname{Span} \{ w_j^{\pm}, \widetilde{w}_j^{\pm}: 
j=0,1,\dots, n\}, \qquad E'' = \operatorname{Span} \{ w_j^{\pm} , 
\widetilde{w}_j^{\pm} : j \geq n+1 \}, $$
and write 
$$\phi = \phi' + \phi'' , \qquad \phi' \in E', \qquad \phi'' \in 
E''.$$
We claim that there exists $c>0$ independent of any 
choices we have made such that 
$|\phi(t,\theta)| \leq c e^{- t}$. We write 
$$\phi' = \sum_{j=0}^n 
(c_j^\pm w_j^\pm + \tilde c_j \widetilde{w}_j^\pm), \qquad 
\phi'' = \sum_{j=n+1}^\infty (c_j^\pm w_j^\pm + \widetilde{c}_j^\pm
\widetilde{w}_j^\pm ).$$ 
By construction $\phi$ decays at some rate as $t \rightarrow \infty$, 
so $c_0^\pm = 0$ and $\widetilde{c}_0^\pm = 0$. Also by 
construction $|\phi(t,\theta)| \leq 1$ for $t\geq 0$, which 
implies $|c_j^\pm|$ and $|\widetilde{c}_j^\pm|$ are all 
bounded independent of all choices for $j=1,\dots, n$. We conclude 
$|\phi'(t,\theta)| \leq c' e^{- t}$ for some 
$c'>0$ which is independent of all choices. We complete the 
proof of the claim by showing $|\phi''| \leq c'' e^{- t}$ 
for some $c''>0$ not depending on any choices. Suppose otherwise, 
then there would exist a sequence $\phi_i'' \in E''$ such that 
$$A_i  =\sup_{t>0} e^{ t} |\phi_i'' (t,\theta)| 
\rightarrow \infty, \qquad \phi_i''(t_i,\theta_i) = 
e^{- t_i} A_i.$$
Define 
$$\widetilde {\phi}_i : [-t_i,\infty) \times \Ss^{n-1} \rightarrow \R, \qquad 
\widetilde{\phi}_i(t,\theta) = \frac{e^{ t_i}}{A_i} 
\phi_i''(t+t_i,\theta) \in E'' .$$
The PDE \eqref{linearization2} is uniformly elliptic on any bounded 
cylinder, so the sequence $\{ t_i\}$ cannot be bounded. Thus $t_i 
\rightarrow \infty$ and (after passing to a subsequence) we obtain 
a limit a $\widetilde \phi \in E''$ on the whole cylinder $\R \times \Ss^{n-1}$. 
Moreover, for $t \geq -t_i$ we have 
$$\widetilde{\phi}_i(t,\theta) = \frac{e^{ t_i}}{A_i} 
\phi_i''(t+t_i,\theta) \leq \frac{e^{t_i}}{A_i} 
e^{- (t+t_i)} A_i = e^{- t},$$
so that in the limit we obtain a global solution $\widetilde{\phi} \in E''$ 
on the whole cylinder $\R \times \Ss^{n-1}$ which satisfies the bound 
$|\widetilde \phi(t,\theta)| \leq e^{- t}$. However, 
this contradicts the definition of $E''$, as all functions in this space 
must grow as fast as $e^{\gamma_{\epsilon,2}|t|}$ either as $t 
\rightarrow \infty$ or $t \rightarrow -\infty$. 
\end{proof} 

\begin{lemma} \label{helper_lemma4} 
Let $v: (0,\infty) \times \Ss^{n-1} \rightarrow (0,\infty)$ 
solve \eqref{paneitz_pde2} and let $\tau_j$ be a sequence 
such that $\tau_j \rightarrow \infty$. Define $v_j(t,\theta) = 
v(t+\tau_j, \theta)$ and 
$$w_j (t,\theta) = v(t+\tau_j,\theta) - v_\epsilon(t+T), \qquad 
\alpha_j = \sup_{0\leq t \leq N T_\epsilon} |w_j(t,\theta)|, \qquad 
\phi_j = \frac{1}{\alpha_j} w_j .$$
Here $v_\epsilon (\cdot + T)$ is the limit extracted from the 
sequence $\{ v_j\}$ using {\it a priori} upper and lower bounds. 
The sequence $\{ \phi_j\}$ converges to a Jacobi field $\phi$ for 
$v_\epsilon ( \cdot + T)$. Moreover, $\phi$ is bounded for 
$t\geq 0$. 
\end {lemma} 

\begin{proof} 
Let $\Delta_{\textrm{cyl}}^2$ denote $\Delta^2$ written in 
cylindrical coordinates, {\it i.e.} 
\begin{equation} \label{cyl_bilap}
\Delta_{\textrm{cyl}}^2 = \frac{\partial^4}{\partial t^4} - \left ( 
\frac{n(n-4)+8}{2} \right ) \frac{\partial^2}{\partial t^2} + 
\frac{n^2(n-4)^2}{16} + \Delta_\theta^2 + 2\Delta_\theta 
\frac{\partial^2}{\partial t^2} - \frac{n(n-4)}{2} \Delta_\theta.
\end{equation} 
Then 
\begin{eqnarray*} 
\Delta_{\textrm{cyl}}^2 w_j & =& \frac{n(n-4)(n^2-4)}{16} 
\left ( v_j^{\frac{n+4}{n-4}} - v_\epsilon ^{\frac{n+4}{n-4}} \right ) \\ 
& = & \frac{n(n-4)(n^2-4)}{16} \left ( \frac{n+4}{n-4} 
\right ) v_\epsilon^{\frac{8}{n-4}} w_j + \mathcal{O}
(\| v_j - v_\epsilon \|^2) 
\end{eqnarray*} 
which proves the sequence $\{ w_j\}$ converges to 
a Jacobi field. It then follows that $\{ \phi_j\}$ also
converges to a Jacobi field $\phi$. 

Next we prove $\phi$ is bounded for $t \geq 0$. Decompose 
$\phi$ into Fourier modes, writing 
$$\phi = \phi'+ \phi'' = a_+w_{\epsilon}^+ + a_-
w_0^- + \phi'', $$
where $\phi''$ is the sum of all the nonzero Fourier modes. We 
may further decompose $\phi'' = \phi''_+ + \phi''_-$, where $\phi''_+$
is the sum of all the Fourier modes which grow exponentially as 
$t \rightarrow + \infty$ and $\phi''_-$ is the sum of Fourier modes 
which decay exponentially as $t \rightarrow +\infty$. Observe that 
$\phi''$ is bounded if and only if $\phi''_+ = 0$. Now fix $\theta \in 
\Ss^{n-1}$ and observe $\partial_\theta \phi''_+$ grows exponentially, 
at least as fast as $e^{\gamma_{\epsilon, 1} t}$, while $\partial \phi''_-$ 
decays exponentially, at least as fast as $e^{- t}$. 
Thus $\phi''$ is bounded in the half-cylinder $t\geq 0$ if and only if 
$\partial_\theta \phi''$ is bounded. 

We can write a similar Fourier decomposition for $\phi_j = \phi_j ' +
\phi_j ''$ for each $j$. Observe that each $\phi_j '$ does not depend at 
all on $\theta$, so that 
$$\partial_\theta \phi_j'' = \partial _\theta \phi_j = \frac{1}{\alpha_j} 
\partial_\theta w_j = \frac{1}{\alpha_j} \partial_\theta v_j.$$
We have just shown in Lemma \ref{helper_lemma3} that 
$\partial_\theta v_j$ converges to a Jacobi field $\psi$ such 
that $|\psi (t,\theta)|\leq c e^{-t}$ for 
$t>0$, where $c$ does not depend on any choices. 
Now let 
$$A_j = \sup \{ |\partial_\theta v_j(t,\theta)| = |\partial_\theta
v(t+\tau_j,\theta)| : 0 \leq t \leq N T_\epsilon \}.$$
If there exist $C_1>0$ and $C_2>0$ such that 
$$C_1 A_j \leq \alpha_j \leq C_2 A_j,$$
{\it i.e.} if $A_j$ and $\alpha_j$ are commensurate, 
then $\partial_\theta \phi'' = \partial_\theta \phi$ is commensurate 
with the Jacobi field $\psi$, and hence decays exponentially. If no uniform 
upper bound $C_2$ exists, then the sequence 
$$\left \{ \frac{1}{\alpha_j} \partial_\theta v_j \right \}$$ 
is unbounded and cannot converge, contradicting the convergence 
$\phi_j \rightarrow \phi$. Additionally, if no lower bound $C_1$ exists
then $\phi_j \rightarrow 0$ uniformly on $[0,N T_\epsilon ] \times 
\Ss^{n-1}$, which is also a contradictions. Thus $\{A_j\}$ and 
$\{ \alpha_j\}$ must be commensurate and so $\phi''$ must 
decay exponentially.  

Lastly we show $a_- = 0$, completing the proof that $\phi$ is 
bounded for $t \geq 0$. We have 
$$v_j (t,\theta) = v_\epsilon (t+T) + \alpha_j \phi + o(\alpha_j),$$
which implies 
\begin{equation} \label{changing_poho}
\mathcal{H}_{\textrm{rad}} (v_j) = \mathcal{H}_{\textrm{rad}} 
(v_\epsilon) + a_-\frac{d}{d\epsilon} \mathcal{H}_{\textrm{rad}}
(v_\epsilon) + o(\alpha_j).
\end{equation}
However, we have already shown $\mathcal{H}_{\textrm{rad}}
(v_j) = \mathcal{H}_{\textrm{rad}} (v_\epsilon)$ for each $j$, 
so \eqref{changing_poho} is only possible if $a_- = 0$.

We have now shown $\phi = a_+ w_0^+ + \phi''$, where $\phi''$ 
decays exponentially. However $w_0^+$ is periodic and hence bounded, 
and so $\phi$ is bounded as well. 
\end{proof} 

\begin{lemma} \label{helper_lemma5}
Let $v: (0,\infty) \times \Ss^{n-1} \rightarrow (0,\infty)$ satisfy 
\eqref{paneitz_pde2} and \eqref{pos_lap_cyl_coords}, let $\tau \geq 0$, 
and let $B>0$. Let 
$$w_\tau(t,\theta) = v_\tau(t,\theta) - v_\epsilon (t+T) = v(t+\tau,\theta) 
- v_\epsilon(t+T)$$
and for $N \in \mathbf{N}$ let 
$$\eta(\tau) = \sup_{0\leq t \leq NT_\epsilon} |w_\tau(t,\theta)|.$$
If $\tau$ is sufficiently large and $\eta(\tau)$ is sufficiently small 
then there exists $s$ such that $|s| \leq B \eta(\tau)$ such that 
\begin{equation} \label{pre_exp_decay1}
\eta (\tau + N T_\epsilon + s) \leq \frac{1}{2} \eta (\tau).
\end{equation} 
\end{lemma} 

\begin{proof} Fix $B >0$ and $N \in \mathbf{N}$. If \eqref{pre_exp_decay1} does 
not hold then there must exist $\tau_j \rightarrow \infty$ such that $\eta_j 
= \eta (\tau_j) \rightarrow 0$ but for any $s$ such that $|s| \leq B\eta_j$ 
it holds $\eta (\tau_j + N T_\epsilon + s) \geq \frac{1}{2} \eta_j$. Let $\phi_j 
= \frac{1}{\eta_j} w_{\tau_j}$. By Lemma \ref{helper_lemma4} the 
sequence $\{ \phi_j\}$ converges uniformly on compact subsets to a 
Jacobi field $\phi$. Moreover, $\phi$ is bounded for $t \geq 0$. By 
hypothesis 
$$0 \leq t \leq N T_\epsilon \Rightarrow |\phi(t,\theta)| \geq \frac{1}{2} $$
so $\phi$ is not identically zero. In the half-cylinder $(0,\infty) \times \Ss^{n-1}$
we expand $\phi = a w_0^+ + \widetilde \phi$ where $\widetilde \phi$ is a 
sum of Fourier modes all of which decay at least like $e^{-\gamma_{\epsilon, 1} t}$. 
The coefficients of the other Fourier modes must all be zero because 
$\phi$ is bounded for $t \geq 0$. 

Next we adjust $v_\epsilon(\cdot + T)$ by the translation parameter 
$s_j = - \eta_j a$. Without loss of generality $B > |a|$, so that $|s_j| \leq 
B \eta_j$. Thus 
$$
w_{\tau_j + s_j} (t,\theta) = v(t+\tau_j - \eta_j a, \theta) - v_\epsilon(t+T)
= w_{\tau_j}(t,\theta) - a \eta_j w_0^+(t) + o(\eta_j), $$
which implies 
$$w_{\tau_j+ s_j} = \eta_j \widetilde \phi + o(\eta_j).$$
The function $\widetilde \phi$ decays at least like $e^{- t}$
and $\sup_{0 \leq t \leq NT_\epsilon} |\widetilde \phi(t,\theta)| \leq 1$, so 
we may choose $N$ sufficiently large such that 
\begin{eqnarray*} 
\sup_{0 \leq t \leq NT_\epsilon} |w_{\tau_j + s_j + N T_\epsilon} (t,\theta)| 
& = & \sup_{N T_\epsilon \leq t \leq 2N T_\epsilon} |w_{\tau_j + s_j} 
(t,\theta)| \\ 
& = & \eta_j \sup_{N T_\epsilon \leq t \leq 2N T_\epsilon} |\widetilde 
\phi(t,\theta) + o(\eta_j) \\ 
& \leq & \frac{1}{4} \eta_j. 
\end{eqnarray*} 
However, this contradicts the hypothesis $\eta(\tau_j + N T_\epsilon + s) 
> \frac{1}{2} \eta(\tau_j)$, completing the proof. 
\end{proof} 

We complete our proof of \eqref{asymp_symm1} by showing we can 
choose $\sigma$ such that $w_\sigma \rightarrow 0$. Again we choose 
$B> 0$ and $N \in \mathbf{N}$. Letting $t> 0$ be sufficiently large we 
may assume $B \eta(0) \leq \frac{1}{2} N T_\epsilon$. Let $\tau_0 = 0$ 
and choose $s_0$ such that \eqref{pre_exp_decay1} holds. 

Next we choose the sequences $\tau_j$, $s_j$, and $\sigma_j$ 
by 
$$\tau_j = \tau_{j-1} + s_{j-1} + N T_\epsilon, \qquad 
\sigma_j = \sum_{i=0}^{j-1} s_i,$$
where $s_j$ satisfies \eqref{pre_exp_decay1} with the choice 
$\tau = \tau_j$. Iterating we see 
$$\eta(\tau_j) \leq 2^{-j} \eta(0) \Rightarrow |s_j| \leq 2^{-j-1} N T_\epsilon 
\Rightarrow \sigma = \lim_{j \rightarrow \infty} \sigma_j = 
\sum_{i=0}^\infty s_j \leq N T_\epsilon. $$
Finally we must show that $\sigma$ is indeed the correct translation 
parameter. For any $t > 0$ write $t = j N T_\epsilon + [t]$ where 
$0 \leq [t]< N T_\epsilon$. Then 
\begin{eqnarray*} 
w_\sigma (t,\theta) & = & v(t+\sigma, \theta) - v_\epsilon (t+T) \\ 
& = & v(t+\sigma_j, \theta) - v_\epsilon (t+T) + v(t+\sigma, \theta) - 
v(t+\sigma_j, \theta) \\ 
& = & w_{\tau_j} ([t],\theta) + \mathcal{O}(2^{-j}). 
\end{eqnarray*} 
Our bound on $\eta(\tau_j)$ then implies
$$v(t+\sigma, \theta) - v_\epsilon (t+T) = 
w_\sigma(t,\theta) \leq C 2^{-j},$$
which  is exactly the exponential decay we claimed. 
\end {proof} 

\section{Refined asymptotics} 
\label{refined_sec}

We finally derive a refined asymptotic expansion of solutions 
of \eqref{paneitz_pde1} with an isolated singularity, essentially 
writing out the next term in the Taylor expansion. 

\begin{thm} 
Let $v:(0,\infty) \times \Ss^{n-1} \rightarrow (0,\infty)$ be a 
smooth solution of \eqref{paneitz_pde2} such that the associated 
metric $g_{ij} = v^{\frac{4}{n-4}} (dt^2 + d\theta^2)$ has positive 
scalar curvature. Then either $\limsup_{t \rightarrow \infty} v(t,\theta)
= 0$, in which case $v$ is a spherical solution, or there exist $\epsilon 
\in [\epsilon_n,1)$, $T\in [0,T_\epsilon)$, $a \in \R^n$, $C>0$ and $\beta>1$ 
such that
\begin{equation} \label{refined_asymp_estimate}
|v(t,\theta) - v_{\epsilon, a} (t+T, \theta)| \leq C e^{-\beta t}.
\end{equation} 
\end{thm} 

\begin{proof} By \eqref{asymp_symm1} there exist $\epsilon$, 
$T$, $C_1$, and $\alpha$ such that 
$$|v(t,\theta) - v_\epsilon (t+T) | \leq C_1 e^{-\alpha t}.$$
In other words 
$$w(t,\theta) = v(t,\theta) - v_\epsilon (t+T) \in 
\mathcal{C}^\infty_{-\alpha} ((0,\infty) \times \Ss^{n-1}).$$
However,  
\begin{eqnarray} 
L_\epsilon (w) & = & \Delta^2_{\textrm{cyl}} (v - v_\epsilon (\cdot + T)) 
- \frac{n(n+4)(n^2-4)}{16} v_\epsilon^{\frac{8}{n-4}} w \\ \nonumber 
& = & \frac{n(n-4)(n^2-4)}{16} \left ( v^{\frac{n+4}{n-4}} - 
v_\epsilon^{\frac{n+4}{n-4}} \right ) 
- \frac{n(n+4)(n^2-4)}{16} v_\epsilon^{\frac{8}{n-4}} w \\ \nonumber 
& = & \mathcal{Q}_{\textrm{cyl}} (w) \in \mathcal{C}^\infty_{-2\alpha}
((0,\infty) \times \Ss^{n-1})
\end{eqnarray} 
where $\Delta^2_{\textrm{cyl}}$ is given in \eqref{cyl_bilap}. We 
combine this with Proposition \ref{lin_soln_asymp_expansion} to 
see $w \in \mathcal{C}^\infty_{-2\alpha} ((0,\infty) \times \Ss^{n-1})$. 
We iterate this finitely many times to obtain $w \in \mathcal{C}^\infty_{-\delta}
((0,\infty) \times \Ss^{n-1})$ for some $\delta \in (1/2, 1)$, and so 
we can apply Proposition \ref{lin_soln_asymp_expansion} once 
more to see 
$$w \in \mathcal{C}^\infty_{-\beta} ((0,\infty) \times \Ss^{n-1}) \oplus 
\operatorname{Span} \{w_j^\pm, j=0,1,2,\dots, n\}, \qquad 
\beta = \min \{ 2\delta, \gamma_{\epsilon,2} \}.$$
However, combining \eqref{trans_del_soln3} and \eqref{low_paneitz_mode2}
we have 
$$v_{\epsilon, a} (t+T,\theta) = v_\epsilon (t+T) + e^{-t-T} \langle \theta, a
\rangle w_1^+(t+T) + \mathcal{O} ( e^{-2t}) ,$$ 
which completes the proof. 
\end{proof}


\end {document}